\documentclass{amsart}
\usepackage{amsmath,amsfonts,epsfig,amscd}
\usepackage[pdftex]{hyperref}

\begin{document}
\bibliographystyle{plain}

\newtheorem{thm}{Theorem}[section]
\newtheorem{lem}[thm]{Lemma}
\newtheorem{prop}[thm]{Proposition}
\newtheorem{cor}[thm]{Corollary}
\newtheorem{conj}[thm]{Conjecture}
\newtheorem{mainlem}[thm]{Main Lemma}
\newtheorem{defn}[thm]{Definition}
\newtheorem{rmk}[thm]{Remark}

\def\square{\hfill${\vcenter{\vbox{\hrule height.4pt \hbox{\vrule width.4pt
height7pt \kern7pt \vrule width.4pt} \hrule height.4pt}}}$}

\newenvironment{pf}{{\it Proof:}\quad}{\square \vskip 12pt}
\newcommand{\dt}{\ensuremath{\text{det}}}
\newcommand{\U}{\ensuremath{\widetilde}}
\newcommand{\Hn}{\ensuremath{\mathbb{H}^3}}
\newcommand{\h}{{\text{hyp}}}
\newcommand{\acyl}{{\text{acyl}}}
\newcommand{\inc}{{\text{inc}}}
\newcommand{\tp}{{\text{top}}}
\newcommand{\V}{\ensuremath{\text{Vol}}}
\newcommand{\Hess}{\ensuremath{{\text{Hess} \ }}}

\title{Hyperbolic convex cores and simplicial volume}
\author{Peter A. Storm}

\date{August 22nd, 2005}

\begin{abstract}
This paper investigates the relationship between the topology of hyperbolizable $3$-manifolds $M$ with
incompressible boundary and the volume of hyperbolic convex cores homotopy equivalent to $M$.  Specifically, it
proves a conjecture of Bonahon stating that the volume of a convex core is at least half the simplicial volume of
the doubled manifold $DM$, and this inequality is sharp.  This paper proves that the inequality is in fact sharp
in every pleating variety of $\text{AH}(M)$.
\end{abstract}

\thanks{This research was partially supported by a National Science Foundation Postdoctoral Fellowship and the Clay Mathematics Institute Liftoff program.}


\maketitle

\section{Introduction}
Let $M$ be a compact oriented $3$-manifold with nonempty incompressible boundary such that the interior of $M$
admits a complete convex cocompact hyperbolic metric.  Define the topological invariant 
$$\mathcal{V}(M) := \inf \{ \text{Vol}(C_N) \, | \, 
  N \text{ hyperbolic \& homeomorphic to int}(M) \},$$ 
where $C_N$ is the convex core of $N$.  Bonahon conjectured that $\mathcal{V}(M)$ equals half the simplicial volume of the doubled manifold $DM$ (see also \cite[Sec.7]{CMT}).  This paper proves Bonahon's conjecture.

The beginning of this paper proves the following Kleinian version of the Besson-Courtois-Gallot theorem
\cite{BCGlong}.
\begin{thm} \label{BCGthm}
If $N$ is a hyperbolic $3$-manifold homotopy equivalent to $M$ then $$\text{Vol}(C_N) \ge \frac{1}{2}
\text{SimpVol}(DM).$$ Moreover, if $\text{Vol}(C_N) = \frac{1}{2} \text{SimpVol}(DM) > 0$ then $M$ is
acylindrical, $N$ is convex cocompact, and $\partial C_N$ is totally geodesic.
\end{thm}
{\noindent}(See Theorems \ref{lower bound thm2} and \ref{equality case thm}.)   Note that $\V (C_N) = 0$ if and only if $N$ is either a Fuchsian group or an extended Fuchsian group.  (An extended Fuchsian group is a degree two extension of a Fuchsian group.)  The main tool used to prove Theorem \ref{BCGthm} is due to Souto \cite{BCS1}, building on work of Besson-Courtois-Gallot \cite{BCGlong}.

Let us postpone an introduction of simplicial volume to Section \ref{simplicial volume section}.  Briefly, if $X$ is a closed oriented $3$-manifold for which Thurston's geometrization conjecture is true, then it is possible to cut $X$ along essential spheres and tori into manifolds admitting one of Thurston's eight geometries.  It is known the simplicial volume of $X$ is the total volume of the hyperbolic pieces of this decomposition of $X$.  Simplicial volume also known as the Gromov norm.  (They differ by a constant.)  

In its second half, this paper proves the inequality of Theorem \ref{BCGthm} is sharp.  In fact, using the work of Bonahon-Otal \cite{BO}, the inequality is proven to be sharp in every pleating variety of the deformation space.\vskip 6pt

\noindent\textbf{Theorem \ref{combinations thm 2}.} \itshape Let $\lambda$ be the bending measured lamination on
$\partial M$ of the boundary of the convex core of a convex cocompact hyperbolic manifold $N$ homeomorphic to the
interior of $M$.  Then for $0 < \varepsilon < 1$ there exists a convex cocompact hyperbolic manifold $N_\varepsilon$ homeomorphic to
the interior of $M$ with bending measured lamination $\varepsilon \lambda$ such that 
$$\V (C_{N_\varepsilon}) \longrightarrow \frac{1}{2} \text{SimpVol}(DM) \quad \text{as} 
     \quad \varepsilon \rightarrow 0. \qquad (\dagger)$$
\normalfont \vskip 6pt 

{\noindent}(The existence of the manifolds $N_\varepsilon$ above follows quickly from the work of Bonahon-Otal \cite{BO}.  The central point of Theorem \ref{combinations thm 2} is $(\dagger)$.)  Together, Theorems \ref{BCGthm} and \ref{combinations thm 2} solve the above conjecture.  Namely,

\begin{cor} \label{BOA2 conj} $\mathcal{V}(M) = \frac{1}{2} \text{SimpVol}(DM).$
In particular, $\mathcal{V}$ is a homotopy invariant of hyperbolizable compact oriented $3$-manifolds with
incompressible non-toroidal boundary.
\end{cor}
{\noindent}(See Corollary \ref{simpvol 3}.)  A pared version of Corollary \ref{BOA2 conj} is also proved.

\vskip 6pt

{\noindent}\emph{Acknowledgements:}  The author thanks Richard Canary for his essential assistance at every stage of this 
research.  Juan Souto's doctoral thesis is a fundamental and essential tool in this paper.  The author thanks him for 
visiting the University of Michigan and sharing some of his ideas.  This research benefited from conversations with 
Ian Agol, Cyril Lecuire, and Caroline Series.  Much of this manuscript was prepared during a stays at the Fourier 
Institute in Grenoble and the Isaac Newton Institute in Cambridge.  The author thanks his hosts for their generous 
hospitality.

\subsection{Outline}

The two main tools used to prove Theorem \ref{BCGthm} are an inequality of Souto \cite[Sec.3,Thm.4]{So}, and a 
previous result of the author \cite[Thm.8.1]{St2}.  For convenience, we restate them here.  For a geometrizable closed 
$3$-manifold $X$, $X_\h$ is the finite volume hyperbolic manifold obtained from the geometric decomposition of $X$.

\begin{thm} \label{Souto's inequality} \cite[Sec.3,Thm.4]{BCS1} \emph{(Souto's inequality)}
If $X$ is a closed geometrizable Riemannian $3$-manifold, then $$h( \U{X})^3 \, \text{Vol}(X) \ge 2^3 \, \V (X_\h)
= 2^3 \, \text{SimpVol}(X),$$ where $h(\U{X})$ is the volume growth entropy of the universal cover of $X$.
\end{thm}

{\noindent}(Souto also proved a related rigidity theorem.  See \cite{BCS1}.)  Let $M$ be a hyperbolizable compact connected oriented $3$-manifold with incompressible boundary.  Let
$\mathcal{I}(M)$ be the set of isometry classes of hyperbolic $3$-manifolds homotopy equivalent to $M$.  A
geometrically finite hyperbolic $3$-manifold $N$ has Fuchsian ends if each component of $N - C_N$ has Fuchsian
holonomy.  It is a consequence of Thurston's Geometrization Theorem and Mostow Rigidity that $M$ is acylindrical
if and only if there exists a unique convex cocompact $N_g \in \mathcal{I}(M)$ with Fuchsian ends (Corollary
\ref{pared acylindrical}).

\begin{thm} \label{my thesis} \cite[Thm.8.1]{St2}
Let $M$ be a compact acylindrical $3$-manifold.  Let $N_g \in \mathcal{I}(M)$ be the unique convex cocompact
manifold with Fuchsian ends.  Then for all $N \in \mathcal{I}(M)$, $$\text{Vol}(C_N) \ge \text{Vol} (C_{N_g}),$$
with equality if and only if $N$ and $N_g$ are isometric.
\end{thm}

We now outline the contents of Sections \ref{pared lower bound}-\ref{combinations section}.  (Section 2 is devoted to preliminaries.)  Define the constant 
$$\mathcal{V} := \frac{1}{2} \text{SimpVol}(DM).$$

\vskip 6pt {\noindent}\emph{Section \ref{pared lower bound}:}  For $N \in \mathcal{I}(M)$ we first prove that $\V(C_N)
\ge \mathcal{V}$.  The rough idea is to double $C_N$ to form $DC_N$, smoothly approximate the
resulting metric space by Riemannian manifolds $A_i^\infty$, and apply Souto's inequality to the pair $A_i^\infty$
and $DM$.  Next, the rigidity statement of Theorem \ref{BCGthm} is proven.  Namely, if 
$N \in \mathcal{I}(M)$ and $\V(C_N) = \mathcal{V} > 0,$ then $N$ is convex cocompact with Fuchsian ends, and $M$ is 
acylindrical.  For this outline suppose that $N$ is convex cocompact.  Assume $M$ is not acylindrical.  Then $DM$ is not hyperbolizable, and consequently $(DM)_\h$ is a \emph{noncompact} collection of finite volume hyperbolic manifolds.  By \cite{L}, there is a sequence of Riemannian metrics $g_n$ on $DM$ converging geometrically to $(DM)_\h$.  In particular, the diameter of $(DM, g_n)$ is going to infinity.  We may use the barycenter method of Besson-Courtois-Gallot and a technique of Souto to produce uniformly Lipschitz surjections $DC_N \longrightarrow (DM, g_n)$.  Since $DC_N$ is compact, this is a contradiction.  Therefore $M$ is acylindrical.  An application of Theorem \ref{my thesis} completes the proof. \vskip 6pt

{\noindent}\emph{Section \ref{gccc section}:}  This section collects results from the literature about how 
continuously convex cores vary under geometric convergence, and the behavior of convex core volumes under geometric 
convergence.  Results due to Bowditch \cite{Bow}, McMullen \cite{Mc3}, and Taylor \cite{Ta} are used.  These facts are 
later used to prove Theorem \ref{combinations thm 2}. \vskip 6pt

{\noindent}\emph{Section \ref{combinations section}:}  This section proves Theorem \ref{combinations thm 2}.  In 
particular, we find a sequence of convex cocompact manifolds $N_n \in \mathcal{I}(M)$ such that $\V (C_{N_n}) 
\rightarrow \frac{1}{2} \text{SimpVol}(DM)$.  For this, the idea is to take manifolds $N_n$ such that the bending 
measure on $\partial C_{N_n}$ is a fixed projective measured lamination with total measure going to zero.  Such 
manifolds were proven to exist by Bonahon-Otal \cite{BO}.  Some work is required to prove that these ``lightly 
pleated'' manifolds converge in some sense to a (disconnected) manifold with Fuchsian ends, and that convex core 
volume tends toward the expected value $\mathcal{V}$. \vskip 6pt

\section{preliminaries}

This section collects the most important definitions for this paper.

\subsection{Pared 3-manifolds \normalfont{\cite[Def.4.8]{M}}} \label{pared 3-manifolds section}
Let $M$ be a compact orientable irreducible $3$-manifold with nonempty boundary.  Assume no component of $M$ is a
$3$-ball.  Let $P \subseteq \partial M$ be a closed subset.  $(M,P)$ is a \emph{pared} $3$-manifold if the
following three conditions hold. \\ (1)  Every component of $P$ is an incompressible torus or a compact annulus.\\
(2)  Every noncyclic abelian subgroup of $\pi_1 (M)$ is conjugate into the fundamental group of a component of
$P$.\\ (3)  Every $\pi_1$-injective cylinder $C : (S^1 \times I, S^1 \times \partial I) \longrightarrow (M,P)$ is
relatively homotopic to a map $\psi$ such that $\psi (S^1 \times I) \subseteq P$.

The submanifold $P \subset \partial M$ is the \emph{paring locus} of $(M,P)$.  By Thurston's Geometrization
Theorem \cite{M}, $(M,P)$ is a pared $3$-manifold if and only if there exists a geometrically finite hyperbolic
structure on the interior of $M$ such that $C_N \cong M - P$.

Many standard topological notions have natural adaptations to the pared category.  For example, a \emph{pared
homotopy equivalence} is (coincidentally) a homotopy equivalence of pairs.  $(M,P)$ has \emph{pared incompressible
boundary} if the inclusion map $\partial M - P \hookrightarrow M$ is $\pi_1$-injective on all components.  An
annulus $C:S^1 \times I \longrightarrow M$ is \emph{essential} in $(M,P)$ if $C$ is $\pi_1$-injective, $C$ is a
map of pairs $$C: (S^1 \times I, S^1 \times \partial I) \longrightarrow (M,\partial M - P),$$ and $C$ is
\emph{not} homotopic rel boundary into $\partial M$ \cite[pg.244]{Th3}.  $(M,P)$ is \emph{pared acylindrical} if
$(M,P)$ has pared incompressible boundary and no essential annuli.  The double $D(M,P)$ of a pared manifold across
its boundary is the compact manifold with toroidal boundary components obtained by gluing two copies of $M$
together along $\overline{\partial M - P}$.  If $(M,P)$ is pared acylindrical then by Thurston's Geometrization
Theorem \cite{M}, $D(M,P)$ admits a complete finite volume hyperbolic metric on its interior.

\subsection{The characteristic submanifold} \label{characteristic submanifold section}

This paper will make extensive use of the characteristic submanifold theory developed by Jaco, Johannson, and Shalen in \cite{JS, Joh, J}.  For a brief discussion of this theory from the point of view of Kleinian groups, the author recommends Morgan's article \cite[pg.88]{M}.  For convenience we here define some of the important terminology.

Let $M$ be a compact oriented irreducible Haken $3$-manifold.  Let $B \subset \partial M$ be an incompressible
subsurface.  There exists a subpair $(\Sigma, S) \subset (M, B)$ with the following properties (see \cite{M, JS,
Joh}):\\ (1)  $\Sigma, S \subset M$ are (possibly disconnected) submanifolds and $S \subset \partial \Sigma \cap
B$. \\ (2)  Each component of $(\Sigma,S)$ is either a $3$-manifold pair of the form $(I-\text{bundle}, \partial
I-\text{subbundle})$ or a Seifert fibered manifold in which the corresponding components of $S$ are foliated by
fibers. \\ (3)  The frontier of $\Sigma$ in $M$ is a collection of essential annuli and tori in $(M,B)$. \\ (4)
No component of $(\Sigma, S)$ is homotopic in $(M, B)$ into another component. \\ (5)  Any essential annulus or
torus in $(M,B)$ is homotopic in $(M, B)$ into $(\Sigma, S)$. \\ Moreover, any pair satisfying these conditions is
isotopic rel $B$ to $(\Sigma, S)$.  By a slight abuse of notation, $(\Sigma, S)$ is called the
\emph{characteristic submanifold} of $(M,B)$.

Now assume $(M,P)$ is a pared manifold with pared incompressible boundary. The characteristic submanifold of
$(M,P)$ as a pared manifold is defined to be the characteristic submanifold of $(M, \overline{\partial M - P})$ as
a pair.  In this case the Seifert fibered components of $\Sigma$ are either solid tori or thickened tori.  A component $L$ of the complement of the characteristic submanifold $(\Sigma,S)$ of $(M,P)$ is either a solid torus, a thickened torus, or $(L, \text{Fr} \Sigma \cap L)$ is pared acylindrical.

The \emph{window} of $(M,P)$ is an $I$-bundle formed as follows.  Begin with the $I$-bundle components of the
characteristic submanifold $\Sigma$ of $(M,P)$.  Now add a regular neighborhood of every essential annulus in the
frontier of $\Sigma$ which is not on the boundary of an $I$-bundle component of $\Sigma$.  Each of these regular
neighborhoods can be viewed as an $I$-bundle with its $\partial I$-subbundle in $S \subset \partial M$.  Some of
these new $I$-bundles may be homotopic rel $\partial M - P$ into other $I$-bundles.  So eliminate any redundant
$I$-bundles.  The resulting collection of $I$-bundles is the \emph{window} of $(M,P)$ \cite{Th5},\cite[Sec.5.3]{CMT}.  The window is an $I$-bundle over a surface in $M$ called the \emph{window base}.

\subsection{Kleinian Deformation theory} \label{deformation theory section}
Let $0< {\mu_3} <1$ be a Margulis constant for hyperbolic $3$-manifolds.  Then for a hyperbolic $3$-manifold $N$,
the  ${\mu_3}$-thin part of $N$ is a disjoint union of bounded Margulis tubes and unbounded cusps \cite{BP}.
After possibly making ${\mu_3}$ smaller, we may also assume that the ends of $\partial C_N$ are totally geodesic
in the $2{\mu_3}$-thin part of $N$ \cite[Lem.6.9]{M}.  Define $N^o$ to be $N$ minus the unbounded components of
its ${\mu_3}$-thin part.  In other words, $N^o$ is the manifold with boundary obtained by removing the cusps from
$N$.  For $\varepsilon < \mu_3$, $N^{\ge \varepsilon}$ (resp. $N^{\le \varepsilon}$) will denote the closed
$\varepsilon$-thick (resp $\varepsilon$-thin) parts of $N$.

Let $M$ be any manifold.  Define the deformation set H$(M)$ as follows.  For an oriented hyperbolic $3$-manifold $N$ and a 
homotopy equivalence $m: M \longrightarrow N$, the pair $(N,m)$ is in the set H$(M)$.  The pairs $(N_1, m_1)$ and 
$(N_2, m_2)$ are equivalent in $H(M)$ if there exists an orientation preserving isometry $\iota: N_1 \longrightarrow 
N_2$ such that $\iota \circ m_1 \sim m_2$.

We will also need a relative version of the above set of hyperbolic structures.  So now let $(M,P)$ be a pared $3$-manifold.  Define the deformation set H$(M,P)$ as follows.  For a hyperbolic $3$-manifold $N$ and a map $m: M \longrightarrow 
N^o$, $(N,m) \in \text{H}(M,P)$ if there exists a union $Q_N$ of components of $\partial N^o$ such that $m: (M,P) 
\longrightarrow (N^o,Q_N)$ is a relative homotopy equivalence. $(N_1, m_1) = (N_2, m_2)$ in H$(M,P)$ if there exists 
an orientation preserving isometry $\iota : N_1 \longrightarrow N_2$ such that $\iota \circ m_1 \sim m_2$.  $(N,m)$ is 
\emph{minimally parabolic} if $Q_N = \partial N^o$.

Let $\mathcal{I}(M,P)$ be the set of (unoriented) isometry classes of hyperbolic $3$-manifolds in $\text{H}(M,P)$.

A hyperbolic manifold $N$ has \emph{Fuchsian ends} if it is geometrically finite and the components of $N - C_N$
have Fuchsian holonomy.  If a connected $N$ has Fuchsian ends, then either $C_N$ or $\partial C_N$ is a totally
geodesic subsurface of $N$.

For precision, it will be convenient to consider hyperbolic manifolds equipped with a framed basepoint.  A framed basepoint $\omega$ of a hyperbolic manifold $N$ is a point $* \in N$ together with an orthonormal basis for the tangent space $T_* N$.  A framed manifold will be a pair $(N,\omega)$.  By fixing for once and for all a framed basepoint $\omega_{\mathbb{H}^3}$ for hyperbolic $3$-space, there is a canonical locally isometric covering map 
$$\iota: (\mathbb{H}^3 , \omega_{\mathbb{H}^3}) \longrightarrow (N, \omega).$$
This covering canonically determines a faithful representation
$$\rho: \pi_1 (N, *) \longrightarrow \text{Isom}^+ (\mathbb{H}^3) \cong \text{PSl}_2 \mathbb{C}.$$
We will use the notation $\pi_1 (N, \omega)$ to denote the subgroup of $\text{PSl}_2\mathbb{C}$ given by the image of $\rho$.

A particularly striking property of pared acylindrical $3$-manifolds is the following corollary of Thurston's
Geometrization Theorem, Mostow rigidity, and Waldhausen's Theorem.
\begin{cor} \label{pared acylindrical} \cite[pg.14]{Th2}
Let $(M,P)$ be a pared acylindrical $3$-manifold.  Then there exists a unique minimally parabolic $N_g \in
\mathcal{I}(M,P)$ with Fuchsian ends.
\end{cor}

\subsection{The bending measured lamination} \label{bending lamination section}

Let $(M,P)$ be a pared $3$-manifold with pared incompressible boundary.  Let $N \in \mathcal{I}(M,P)$ be a
minimally parabolic geometrically finite hyperbolic $3$-manifold equipped with a homeomorphism $$h: M-P
\longrightarrow C_N.$$ The geometry of the boundary of the convex core $C_N$ is a pleated surface, and it is
described by a measured geodesic lamination on $\partial C_N$ \cite{EM}.  Pull back this measured lamination via $h$ to a
measured lamination on $\partial M - P$.  This measured lamination $\beta$ on $\partial M - P$ is the
\emph{bending measured lamination} of the hyperbolic manifold $N$ (homeomorphically) marked by $h$.

In Section \ref{combinations section} we will need the following theorem, whose proof follows from the work of Bonahon \cite{Bon}, Bonahon-Otal \cite{BO}, Kerckhoff-Hodgson \cite{HK}, and Lecuire \cite{Lc}.

\begin{thm} \label{BO thm}
Retain the notation of the previous paragraph.  Pick $0 < \varepsilon < 1$ and consider the measured lamination
$\varepsilon \beta$ on $\partial M - P$.  There exists a minimally parabolic geometrically finite $N_\varepsilon
\in \mathcal{I}(M,P)$ equipped with a homeomorphism $$h_\varepsilon: M-P \longrightarrow C_{N_\varepsilon}$$ such
that:\\
(1)  $\varepsilon \beta$ is the bending measured lamination of $N_\varepsilon$ marked by $h_\varepsilon$, and \\
(2)  The volumes $\V (C_{N_\varepsilon})$ are non-increasing as $\varepsilon \rightarrow 0$.
\end{thm}

Let $p \subset P \subset \partial M$ be an embedded $1$-manifold which is a deformation retract of the annular
components of $P$.  Make $p$ into a measured lamination on $\partial M$ by assigning weight $\pi$ to each
component of $p$.  In order to apply \cite{BO,Lc}, we need
\begin{lem} \label{BO prop}
Let $\lambda$ be a measured lamination on $\partial M - P$.  Consider the following three conditions on $\lambda$:
\\ 
{\noindent}(a)  Each closed leaf of $\lambda$ has weight less than or equal to $\pi$. \\ 
{\noindent}(b)  There exists an $\eta > 0$ such that if $A$ is an essential annulus of $(M,P)$ then $i(\partial A, \lambda) \ge \eta$.\\
{\noindent}(c)  If $D$ is an essential disk of $M$ ($M$ may have compressible boundary!), then $i(\lambda \cup p,
\partial D) > 2\pi$. \vskip 3pt 
If $\lambda$ satisfies conditions (a), (b), and (c) and $0 < \varepsilon < 1$, then $\varepsilon \lambda$ satisfies conditions (a), (b), and (c).
\end{lem}

In the proof of Lemma \ref{BO prop} we will need

\begin{lem} \label{enough intersections}
If $D$ is an essential disk of $M$ then $\partial D \subset \partial M$ intersects intersects $p$ in at least two
points.
\end{lem}
\emph{Proof of Lemma \ref{enough intersections}:} Suppose $\partial D$ intersects $p$ in a single point.  Let $\mathcal{N}p$ (resp. $\mathcal{N}D$) be a regular neighborhood of $p$ (resp. $D$) in $M$.  The frontier of $\mathcal{N}p \cup \mathcal{N}D$ is a compressing disk of $M$ which does not intersect $p$.  Since $(M,P)$ is pared incompressible, this is a contradiction. \square \vskip 6pt

\emph{Proof of Lemma \ref{BO prop}:} Suppose $\lambda$ is a measured lamination of $\partial M -P$ satisfying
conditions (a), (b), and (c).  Then $\varepsilon \lambda$ immediately satisfies condition (a).  Condition (b) is
satisfied for the constant $\varepsilon \cdot \eta$.  Finally, let $D$ be an essential disk of $M$.  By Lemma
\ref{enough intersections}, $i(p, \partial D) \ge 2 \pi$.  Condition (c) applied to $\lambda$ implies that either
$i(p, \partial D) > 2 \pi$ or $i(\lambda, \partial D)>0$.  Therefore 
$$i((\varepsilon \lambda) \cup p, \partial D) = \varepsilon \, i(\lambda, \partial D) + i(p, \partial D) > 2\pi.$$ 
\square \vskip 6pt

\emph{Proof of Thm. \ref{BO thm}:} The proof is broken into two cases.

In the first case let us assume that the support of the measured lamination $\beta$ is a union of simple closed curves.  Then conclusion (1) of Theorem \ref{BO thm} follows from Lemma \ref{BO prop} and \cite[Thm.2]{BO}.  It follows from \cite[Thm.3]{BO} that the manifolds $N_\varepsilon$ are the unique manifolds satisfying conclusion (1).  By \cite[Thm.4.7]{HK}, we may choose the $N_\varepsilon$ and framings $\omega_\varepsilon$ such that the path
\begin{eqnarray*} 
\eta: (0,1)  & \longrightarrow &  \text{Hom}(\pi_1 (M), \text{PSl}_2 \mathbb{C})\\
 \varepsilon & \longmapsto &
   \left\{ \rho_{\varepsilon } : \pi_1 (M) \longrightarrow \pi_1 (N_\varepsilon, \omega_\varepsilon) \right\}
\end{eqnarray*}
is smooth.  Since the path $\eta$ is smooth, a version of the Schl{\"a}fli formula due to Bonahon 
\cite[Cor.2]{Bon} implies that the volume $\V (C_{N_\varepsilon})$ is strictly decreasing as $\varepsilon 
\rightarrow 0$. 

The second case is the general case.  The main tool is \cite[Thm.A]{Lc}.  Using Lemma \ref{BO prop}, we may apply \cite[Thm.A]{Lc} (and its proof in \cite[Sec.5.1]{Lc}) to the lamination $\varepsilon \beta \cup p$.  From this we may conclude that there exists a sequence of measured laminations $\beta^j$ converging to $\beta$ in the space of measured laminations on $\partial M -P$ such that the following are true for all $0 < \varepsilon \le 1$: \\
$\bullet \ \ $  The support of $\beta^j$ is a union of simple closed curves.  \\
$\bullet \ \ $  There exists a sequence of geometrically finite hyperbolic manifolds $N^j_\varepsilon$ with bending measured lamination $\varepsilon \beta^j$ and a homeomorphism $h^j_\varepsilon: M - P \longrightarrow C_{N_\varepsilon^j}$. \\
$\bullet \ \ $  For each $\varepsilon > 0$, any subsequence of the sequence $\{ (N^j_\varepsilon, h^j_\varepsilon)\}_j$ has a further subsequence which converges algebraically to a manifold satisfying conclusion (1) of Theorem \ref{BO thm}.

By a diagonalization argument we can pass to a subsequence $\{ j_k \}$ such that for all rational $\varepsilon \in (0,1)$ the sequence $\{ (N^{j_k}_\varepsilon, h^{j_k}_\varepsilon) \}_k$ converges algebraically to a manifold $( N_\varepsilon, h_\varepsilon)$ satisfying conclusion (1) of Theorem \ref{BO thm}.  For each irrational $\varepsilon \in (0,1)$ define $(N_\varepsilon, h_\varepsilon)$ to be an accumulation point of the sequence $\{ (N^{j_k}_\varepsilon, h^{j_k}_\varepsilon) \}_k$.  To simplify the notation, let us re-index the subsequence $\{ j_k \}$ by the natural numbers, which we will again denote by simply $\{ j \}$.

The geometrically finite manifolds $(N_\varepsilon, h_\varepsilon)$ and $(N^j_\varepsilon, h^j_\varepsilon)$ are minimally parabolic in $\text{H}(M,P)$.  Therefore any algebraically convergent subsequence of $(N^j_\varepsilon, h^j_\varepsilon)$ satisfies the property: $h_\varepsilon (\gamma) \in \text{PSl}_2 \mathbb{C}$ is parabolic if and only if $h_\varepsilon^j (\gamma)$ is parabolic for all $j$.  Such an algebraically convergent sequence is said to be type-preserving.  In the setting of geometrically finite manifolds, it is known that type-preserving algebraically convergence sequences converge strongly.  (A more general theorem is proven by Anderson-Canary in \cite{AC}[Thm.3.1].  Without explicit mention, this issue of strong convergence is carefully studied in \cite{BO}[Ch.2].) 

By the first case in this proof, the manifolds $N^j_\varepsilon$ are uniquely determined and 
$$\varepsilon < \varepsilon' \quad \text{implies} \quad
\V (C_{N^j_\varepsilon}) < \V (C_{N^j_{\varepsilon'}}).$$  
Pick a rational number $t \in (\varepsilon, \varepsilon')$.  By the previous paragraph, the manifold $( N_\varepsilon, h_\varepsilon)$ is a strong accumulation point of the sequence $\{ (N^j_\varepsilon, h^j_\varepsilon) \}$, and for each $j$ we have the inequality
$$\V (C_{N^j_\varepsilon}) < \V (C_{N^j_t}).$$  Convex core volume is continuous in the strong topology \cite{Ta}.  Therefore
$$\V (C_{N_\varepsilon}) \le \V(C_{N_t}).$$
Applying an identical argument to $\varepsilon'$ shows that
$$\V (C_{N_\varepsilon}) \le \V(C_{N_{\varepsilon'}}).$$  This proves part (2) of Theorem \ref{BO thm}.  (The extra argument in this case is necessary because the author does not know how to find a smooth path of deformations analogous to $\eta$.  Without such a smooth path, Bonahon's version of the Schl{\"a}fli formula cannot be applied.) \square \vskip 6pt

\subsection{Lightly bent curves in hyperbolic $3$-manifolds}

In Section \ref{combinations section} it will be necessary to compare the length of a lightly bent essential closed curve on the boundary of the convex core with the length of the curve's geodesic representative.  To do so we use the following estimate due to Lecuire \cite{Lc} (proved also by Series).

\begin{lem} \label{Lecuire lemma} \cite[Lem.A.1]{Lc}
Let $N$ be a geometrically finite hyperbolic $3$-manifold.  Let $\beta$ be the bending measured lamination on the boundary of the convex core $C_N$.  Let $c \subset \partial C_N$ be a closed curve which is a geodesic in the intrinsic hyperbolic metric on $\partial C_N$.  If it exists, let $c^*$ be the closed geodesic of $N$ freely homotopic to $c$.  For $\varepsilon < \pi/2$ there exist constants $K_\varepsilon, A_\varepsilon$ such that: if $i (c, \beta) \le \varepsilon$ then
$$\ell_M (c) \le K_\varepsilon \cdot ( \ell_M (c^*) + A_\varepsilon ),$$
where $K_\varepsilon \rightarrow 1$ and $A_\varepsilon \rightarrow 0$ as $\varepsilon \rightarrow 0$, and $\ell_M (c^*) := 0$ when $c$ is not homotopic to a closed geodesic.  
\end{lem}

\subsection{Simplicial volume} \label{simplicial volume section}

For detailed definitions and an introduction to simplicial volume in the context of $3$-manifolds, the author recommends \cite[Ch.6]{Th}.    

Let $X$ be a closed oriented $3$-manifold.  The ostensibly simple definition of simplicial volume is the following.
$$\text{SimpVol}(X)= V_3 \ \inf \left\{ \sum_i | c_i | \right\},$$
where the infimum is taken over all singular $3$-chains $\sum_i c_i$ realizing the fundamental class of $X$ in singular homology, and $V_3 \approx 1.01$ is the supremal volume of a hyperbolic simplex in $\mathbb{H}^3$.  The simplicial volume divided by $V_3$ is known as the Gromov norm.  There is a similar definition of simplicial volume for compact $3$-manifolds whose boundary consists of tori (see \cite[Ch.6]{Th}).  These definitions will not be used directly in this paper.  It will be more useful to use the relation between simplicial volume and hyperbolic volume which we now describe.

Let $X$ now be a compact oriented irreducible (possibly disconnected) $3$-manifold with (possibly empty) boundary equal to a collection of incompressible tori.  Assume $X$ is \emph{geometrizable}, meaning $X$ can be cut along
embedded essential tori into pieces which admit a geometric structure locally modelled on one of the eight Thurston
geometries.  Let $X_\h$ be the finite volume (possibly disconnected) complete hyperbolic manifold obtained from
the hyperbolic pieces of the geometric decomposition of $X$.  In this situation there is the following useful theorem of Thurston/Soma relating simplicial volume to hyperbolic volume.

\begin{thm} \label{simpvol theorem} \cite{So}\cite[Ch.6]{Th} Let $X$ be as above.  The \emph{simplicial volume} of $X$ 
is the volume of $X_\h$, i.e. $$\text{SimpVol}(X) = \V (X_\h).$$ In particular, if a (possibly disconnected) manifold 
$X'$ is obtained by cutting $X$ along an essential embedded torus, then $$\text{SimpVol}(X) = \text{SimpVol}(X').$$ 
\end{thm}

Let $(M,P)$ be a compact oriented irreducible (possibly disconnected) pared $3$-manifold with pared incompressible
boundary.  Let $\Sigma$ be the characteristic submanifold of $(M,P)$.  Let $(M_\acyl, Q)$ be the compact (not
necessarily connected) pared acylindrical $3$-manifold obtained from the disjoint union of the pared acylindrical
pieces of $M - \Sigma$.    By Corollary \ref{pared acylindrical}, there exists a geometrically finite minimally
parabolic hyperbolic $3$-manifold $N_g \in \mathcal{I}(M_\acyl, Q)$ homeomorphic to $\text{int}(M_\acyl)$ such
that the boundary of the convex core $C_{N_g}$ is totally geodesic.  (Keep in mind that neither $M_\acyl$ nor
$N_g$ is necessarily connected.)

\begin{prop} \label{simpvol 2}
$(D(M,P))_\h$ is isometric to $DC_{N_g}$.  In particular, the simplicial volume of the double $D(M,P)$ is twice
the volume of $C_{N_g}$.
\end{prop}
\begin{pf}
One can verify that the inclusion map $D(M_\acyl,Q) \hookrightarrow D(M,P)$ is $\pi_1$-injective.  Almost by
definition, the complement $D(M,P) - D(M_\acyl,Q)$ can be cut along embedded essential tori into Seifert fiber
spaces.  Therefore $(D(M,P))_\h$ is homeomorphic to the interior of $D(M_\acyl, Q)$, which is in turn homeomorphic
to $DC_{N_g}$.  The proposition now follows from Mostow rigidity and Theorem \ref{simpvol theorem}.
\end{pf}

\begin{cor} \label{simpvol 3}
$\text{SimpVol}(D(M,P))$ is a pared homotopy invariant in the category of compact oriented irreducible pared
manifolds with pared incompressibly boundary.
\end{cor}
\begin{pf}
By Johannson's deformation theorem \cite{Joh}, the pared homeomorphism type of the complement of the characteristic
submanifold is a pared homotopy invariant in the above category.
\end{pf}

Roughly speaking, drilling curves out of a $3$-manifold cannot decrease its simplicial volume.  We shall need a
precise formulation of this.  The following inequality was proved by Agol, completing a sketch in \cite[Ch.6]{Th}.

\begin{prop} \cite{A} \label{simpvol goes up with drilling 2}
Let $X$ be a compact oriented irreducible geometrizable $3$-manifold with (possibly empty) boundary consisting
of incompressible tori.  If $\gamma \subset X$ is a compact embedded $1$-manifold with a regular neighborhood
$\mathcal{N}\gamma$, and $X - \mathcal{N}\gamma$ is hyperbolizable then $$\text{SimpVol}(X - \mathcal{N}\gamma) >
\text{SimpVol}(X).$$
\end{prop}

{\noindent}In \cite{A}, Proposition \ref{simpvol goes up with drilling 2} is stated for $X$ without boundary, but
the proof goes through unchanged for manifolds with toroidal boundary components.  In the setting of pared
manifolds, we have

\begin{cor} \label{simpvol goes up with drilling}
If $(M,P)$ and $(M,Q)$ are pared $3$-manifolds with pared incompressible boundary, and $P \subseteq
Q$ is a collection of connected components, then $$\text{SimpVol}(D(M,Q)) \ge \text{SimpVol}(D(M,P)).$$
\end{cor}
\begin{pf}
A manifold homeomorphic to $D(M,P)$ can be produced by appropriately gluing solid tori onto the boundary of
$D(M,Q)$.  The corollary then follows from \cite[Prop.6.5.2]{Th}.
\end{pf}

\subsection{Volume growth entropy}
Let $X$ be a geodesic metric space of Hausdorff dimension $n$, $\widetilde{X}$ be the universal cover of $X$, and
$\mathcal{H}^n$ be $n$-dimensional Hausdorff measure.  The \textit{volume growth entropy} of $\U{X}$ is the number
$$h(\U{X}) :=  \limsup_{R \rightarrow \infty}
   \frac{1}{R} \log \mathcal{H}^n ( B_{\widetilde{X}} (x,R)),$$
where $x$ is any point in $\widetilde{X}$, and the ball $B_{\widetilde{X}} (x,R)$ is in $\widetilde{X}$. \vskip
3pt {\noindent}The following inequality will be used in Section \ref{pared lower bound}.

\begin{thm}
\label{Perel'man} \cite[pg.40]{BGP} If $X$ is an Alexandrov space with curvature bounded below by $-1$ and
Hausdorff dimension $n$, then the volume growth entropy of $\U{X}$ is less than or equal to the volume growth
entropy of $\mathbb{H}^n$.  In other words $$h(\U{X}) \le h(\mathbb{H}^n) = n-1.$$
\end{thm}

\subsection{Smoothing and geometrically doubling the convex core} \label{convex core smoothings}

In Section \ref{pared lower bound} it will be necessary to consider smooth approximations
of the convex core.  Here we establish some convenient notation and facts.  The geometric structure theory presented in Morgan's article \cite[Sec.6]{M} (in particular \cite[Lem6.8]{M}) will be used repeatedly.

Suppose $N$ is a geometrically finite hyperbolic $3$-manifold. Then $C_N$ has finite volume. For $\varepsilon > 0$, let $\mathcal{N}_\varepsilon C_N$ be an $\varepsilon$-neighborhood of $C_N$.  By \cite[Lem.1.3.6]{EM}, the boundary of $\mathcal{N}_\varepsilon C_N$ is $\mathcal{C}^{1,1}$-smooth.  

For $\delta >0$, let $C_{N,\delta} \subset N$ be a closed convex submanifold of $N$ such that 
$$C_N \subset \text{int}(C_{N, \delta}) \subset C_{N,\delta} \subset \mathcal{N}_\delta C_N,$$
$\V (C_{N,\delta}) < \V(C_N) + \delta$, 
$\partial C_{N, \delta}$ is $\mathcal{C}^{1,1}$-smooth, 
and the unbounded components of $\partial C_{N, \delta} \cap N^{\le \mu_3}$ are totally geodesic.  
One can find such a submanifold by beginning with $\mathcal{N}_\varepsilon C_N$ (for some $\varepsilon < \delta$) and smoothly tapering the ends of $\partial \mathcal{N}_\varepsilon C_N$ so that the ends of the boundary of the tapered manifold are totally geodesic in the ${\mu_3}$-thin part of $N$.  (There are many possible submanifolds satisfying these
conditions.  It will not matter which is chosen.)

We may metrically double $C_{N,\delta}$ across its boundary to obtain the path metric space $DC_{N,\delta}$. Since
$\partial C_{N, \delta}$ is $\mathcal{C}^{1,1}$, the doubled manifold $DC_{N,\delta}$ is a $\mathcal{C}^{1,1}$-manifold without boundary.  The Riemannian metric on $N$ induces a Riemannian metric tensor on $DC_{N, \delta}$, but on a $\mathcal{C}^{1,1}$ manifold it only makes sense to ask whether or not the resulting metric tensor is Lipschitz.  
The induced metric tensor is Lipschitz on each half of the doubled convex core, and glues continuously across the boundary. (Verify this using coordinates in the form of a $\mathcal{C}^{1,1}$ product neighborhood of $\partial C_{N, \delta}$.)  The resulting metric tensor on $DC_{N,\delta}$ is therefore Lipschitz.  Moreover, outside a
compact set $DC_{N, \delta}$ is a hyperbolic rank 2 cusp.  This follows because, by definition, the boundary of $C_{N, \delta}$ is totally geodesic outside a compact set.  (See the previous paragraph for the definition of $C_{N, \delta}$.)  We note that
$DC_{N,\delta}$ is an Alexandrov space with curvature bounded below by $-1$ \cite[Lem.5.4]{St1} (see also
\cite[Thm.5.2]{P}).

Since we are working in dimension three, the $\mathcal{C}^{1,1}$-structure on $DC_{N, \delta}$ contains a unique $\mathcal{C}^\infty$ structure.  So $DC_{N, \delta}$ admits the structure of a smooth manifold equipped with a Lipschitz Riemannian metric.  It is convenient to approximate the Lipschitz metric tensor on $DC_{N, \delta}$ by a smooth Riemannian metric.  A standard smooth approximation argument yields the following lemma, whose proof we omit.

\begin{lem} \label{smooth metric lemma}
Let $(M,g)$ be a smooth manifold equipped with a Lipschitz Riemannian metric tensor.  Pick $L>1$. There
exists a smooth Riemannian metric $h$ on $M$ such that the identity map $$(M, g) \longrightarrow (M, h)$$ is
$L$-Lipschitz.

\end{lem}

\section{Applying the barycenter method} \label{pared lower bound}
Let $(M, P)$ be a connected pared $3$-manifold with pared incompressible boundary.  To avoid trivial cases, we assume that $\text{SimpVol}(D(M,P)) > 0$.  Let $\Sigma$ be the
characteristic submanifold of $(M,P)$.  Let $(M_\acyl, Q)$ be the compact (not necessarily connected) pared
acylindrical $3$-manifold obtained from the disjoint union of the pared acylindrical pieces of $M - \Sigma$.    By
Corollary \ref{pared acylindrical}, there exists a geometrically finite minimally parabolic hyperbolic
$3$-manifold $N_g \in \mathcal{I}(M_\acyl, Q)$ homeomorphic to $\text{int}(M_\acyl)$ such that the boundary of the
convex core $C_{N_g}$ is totally geodesic.  (Keep in mind that neither $M_\acyl$ nor $N_g$ is necessarily
connected.)  Recall that by Proposition \ref{simpvol 2}, $\frac{1}{2} \text{SimpVol}(D(M,P)) = \text{Vol}(C_{N_g})$.  The inequality $\text{SimpVol}(D(M,P))> 0$ implies that $(M_\acyl,Q)$ is not empty.

The goal of this is section is to prove the following two theorems. 
\begin{thm}\label{lower bound thm2}
Let $(M,P)$ and $(N_g, m_g)$ be as above.  If $N \in \mathcal{I} (M,P)$ then $$\text{Vol}(C_N) \ge \V (C_{N_g}) =
\frac{1}{2} \text{SimpVol}(D(M,P)).$$
\end{thm}

\begin{thm}\label{equality case thm}
Let $(M,P)$ and $(N_g, m_g)$ be as above.  If $N \in \mathcal{I}(M,P)$ and $$\text{Vol}(C_N) = \V (C_{N_g}) =
\frac{1}{2} \text{SimpVol}(D(M,P))$$ then $(M,P)$ is pared acylindrical and $N$ is minimally parabolic.
\end{thm}

{\noindent}At the end of this section, Theorem \ref{BCGthm} will be proven using Theorems \ref{lower bound thm2} and \ref{equality case thm}.

Recall there is a natural way to double the pared manifold $(M,P)$ across $\overline{\partial M - P} \subset
\partial M$ to obtain a compact manifold $D(M,P)$ with boundary a disjoint collection of tori (see Section
\ref{pared 3-manifolds section}).  Also double the pared manifold $(M_\acyl,Q)$ across its boundary to obtain
$D(M_\acyl, Q)$.  The doubled manifold $D(M_\acyl, Q)$ is hyperbolizable (see Section \ref{pared 3-manifolds
section}).  Double the compact manifold $D(M,P)$ across its boundary to obtain the closed manifold $DD(M,P)$.  (To
handle all cases simultaneously, we adopt the convention that the double $DX$ of a closed manifold $X$ is simply
two disjoint copies of $X$.)  From Proposition \ref{simpvol 2} it follows that $(DD(M,P))_\h$ is homeomorphic to
two disjoint copies of $\text{int}(D(M_\acyl, Q))$.

Pick a hyperbolic $3$-manifold $N$ in $\mathcal{I} (M,P)$.  Without a loss of generality we may assume that $N$
is geometrically finite.  \emph{Let us temporarily assume that $N$ is minimally parabolic.}  This assumption will be removed later.  Since $\text{SimpVol}(D(M,P))$ is a pared homotopy invariant in the category of pared manifolds with
pared incompressible boundary (Corollary \ref{simpvol 3}), we may assume there exists a homeomorphism 
$$g:(M,P) \longrightarrow (C_N \cap N^{\ge \zeta}, C_N \cap \partial N^{\ge \zeta})$$ for some small $\zeta < \mu_3$ (see Section \ref{deformation theory section}).  The constant $\zeta$ is to remain fixed for the remainder of this section.
  
Metrically double $C_{N,\delta}$ across its boundary to obtain the path metric space $DC_{N, \delta}$.  
(The following facts are discussed in Section \ref{convex core smoothings}.)  
$DC_{N, \delta}$ is a smooth manifold equipped with a Lipschitz
metric tensor.  Outside of a compact set $DC_{N, \delta}$ is a rank $2$ hyperbolic cusp.  
There is a homeomorphism $D(M - P) \longrightarrow DC_{N, \delta}$.  To control
the volume growth entropy of the universal cover of $DC_{N,\delta}$, we need the fact that $DC_{N,\delta}$ is an
Alexandrov space with curvature bounded below by $-1$.

To apply Souto's machinery, we must be working with compact spaces.  If $DC_{N,\delta}$ is not compact, the
following lemma will allow us to make a compact space out of $DC_{N,\delta}$.  Its proof is similar to the proof
of \cite[Prop.2.3]{L} (see also \cite{Bes1} and \cite[Prop.8]{BCS1}).
\begin{lem} \label{gluing cusps lemma}
Let $(T \times [0,\infty), g)$ be a rank $2$ hyperbolic cusp, where $T$ is a $2$-torus.  Pick $\varepsilon > 0$.
There exists an $R_\varepsilon > 0$ and a Riemannian manifold $(T \times [-1,1], g_\varepsilon)$ such that:\\ (1)
The sectional curvature of $g_\varepsilon$ is between $-1-\varepsilon$ and $0$.\\ (2)  The diffeomorphic
involution of $T \times [-1, 1]$ taking $(x,t)$ to $(x,-t)$ is a $g_\varepsilon$-isometry.\\ (3)  There is an
isometric homeomorphism $(T \times [0, R_\varepsilon],g) \longrightarrow (T \times [-1, -1/2], g_\varepsilon)$.\\
(4)  $\V (T \times [-1/2, 1/2], g_\varepsilon) < \varepsilon$. \\ Moreover, $\lim_{\varepsilon \rightarrow 0}
R_\varepsilon = \infty$.
\end{lem}
{\noindent}Recall that $\partial C_{N,\delta} - N^{\ge \zeta}$ is a (possibly disconnected) totally geodesic surface (Section
\ref{convex core smoothings}).  For any $\eta < \zeta$ (where $\zeta < \mu_3$ is the constant fixed above), let
$DC_{N,\delta}^{\ge \eta}$ indicate the compact manifold with horospherical boundary obtained by ``doubling" the
$\eta$-thick part $C_{N,\delta} \cap N^{\ge \eta}$ of $C_{N,\delta}$.  Outside a compact set, $DC_{N,\delta}$ is a smooth hyperbolic rank $2$ cusp.  Therefore by using the above lemma there exists a sequence of path metric spaces
$A_i$ homeomorphic to $DD(M,P)$ such that:\\ 
(a) $A_i$ is an Alexandrov space with curvature bounded below by $-1
- \frac{1}{i}$.\\ 
(b) $A_i$ is a $\mathcal{C}^{1,1}$-manifold with a Lipschitz Riemannian metric tensor.\\ 
(c) There are two isometric embeddings $\phi^k_i: DC_{N,1/i}^{\ge 1/i} \longrightarrow A_i$ (for $k=1,2$) with disjoint images. \\ 
(d) $\V(A_i) \longrightarrow 2 \V (DC_N) = 4 \V (C_N)$. \\

{\noindent}Note that if $C_N$ is compact, then $A_i$ is simply two disjoint copies of $DC_{N,1/i}$.

By Theorem \ref{Perel'man} and property (a) above, $\limsup{h(\U{A_i})} \le h(\Hn) =2$.  As a final approximation,
let $A^\infty_i$ be a $\mathcal{C}^\infty$-smooth Riemannian manifold equipped with $(1+ 1/i)$-bilipschitz
homeomorphism $\psi_i: A_i \rightarrow A^\infty_i$.  (Such a space $A^\infty_i$ exists by Lemma \ref{smooth metric lemma}.)
Then also $$\limsup{h(\U{A^\infty_i})} \le h(\Hn) =2.$$  

We are now prepared to quickly prove Theorem \ref{lower bound thm2}.

\vskip 6pt
\emph{Proof of Thm.\ref{lower bound thm2}:} By Souto's theorem \cite[Sec.3,Thm.4]{BCS1},
$$h(\U{A^\infty_i})^3 \cdot \V (A^\infty_i) \ge 2^3 \, \V ((DD(M,P))_\h) =
             2^3 \cdot 2 \, \V (DC_{N_g}).$$
Taking $i \rightarrow \infty$ and dividing both sides of the inequality by $2^5$ yields $$\V(C_N) \ge
\V(C_{N_g}).$$  This proves Theorem \ref{lower bound thm2} under the additional assumption that $N$ is minimally parabolic.  It remains to remove this assumption.

Recall the homeomorphism
$$g:(M,P) \longrightarrow (C_N \cap N^{\ge \zeta}, C_N \cap \partial N^{\ge \zeta}).$$
Using this homeomorphism, we see that by enlarging the paring locus $R \supset P$, $N$ becomes minimally parabolic in the deformation set $\mathcal{I}(M,R)$ (see Section \ref{deformation theory section}).  Enlarging the paring locus of $(M,P)$ does not decrease the simplicial volume of the double $D(M,P)$ (Corollary \ref{simpvol goes up with drilling}).  This proves the desired inequality in the general case. \square \vskip 6pt

We now prove Theorem \ref{equality case thm}. \vskip 6pt

\emph{Proof of Thm.\ref{equality case thm}:}
Assume that $\V (C_N) = \V(C_{N_g})$.

To begin, we need the following weak form of Ahlfors regularity for $A_i^\infty$, which follows from the Bishop-Gromov
inequality for Alexandrov spaces \cite[Thm.10.6.6]{BBI}.  Recall that for $k \in \{1,2 \}$ and $i \in \mathbb{N}$ we have bilipschitz embeddings
$$\psi \circ \phi_i^k : DC_{N, 1/i}^{\ge 1 / i} \longrightarrow A_i^\infty,$$
where the bilipschitz constant is decreasing to $1$ as $i \rightarrow \infty$.
\begin{lem} \label{Ahlfors reg 2}
There exists a constant $a>0$ such that for all $r<1, i \gg 0, k \in \{1,2\},$ and $x \in DC_{N,1/i}^{\ge \zeta}$
we have $$\V (B_{A_i^\infty} (\psi_i \circ \phi^k_i (x), r) \ge a r^3.$$
\end{lem}
\begin{pf}
Using the bilipschitz embeddings $\psi \circ \phi_i^k: DC_{N,1/i}^{\ge \zeta} \longrightarrow A_i^\infty$, it suffices to prove the lemma
for $DC_{N,1/i}^{\ge \zeta}$.  Since $DC_{N,1/i}^{\ge \zeta}$ is compact, the volume of a unit radius ball in
$DC_{N,1/i}^{\ge \zeta}$ has a positive lower bound which is uniform for $i \gg 0$.  Using this, one may apply the
Bishop-Gromov inequality \cite[Thm.10.6.6]{BBI} to obtain the desired result.  (Alternatively, the Bishop-Gromov inequality can be avoided by directly considering the geometry of $C_{N,1/i}$.)
\end{pf}

We now need the following version of a theorem of Leeb \cite[Prop.2.3,Prop.2.6]{L} (see also \cite[Prop.8]{BCS1}
and \cite{Bes1}).  Recall that $D(M,P)$ is the compact manifold with (possibly empty) toroidal boundary obtained
by doubling $M$ along $\overline{\partial M - P}$, and $DD(M,P)$ is the closed manifold obtained by doubling
$D(M,P)$.
\begin{thm} \label{Leeb's result} \cite{L} (see also \cite[Prop.8]{BCS1})
$DD(M,P)$ admits a family of metrics $\rho_i \ (i \in \mathbb{N})$ with sectional curvature between $-1 -
\frac{1}{i}$ and $0$ such that: 

{\noindent}(1)  There is an open $\pi_1$-injective subset $V \subseteq DD(M,P)$ homeomorphic to $(DD(M,P))_\h$ such that: for all $i$ there is an isometric embedding 
$$(\mathcal{N}_i V, \rho_i) \longrightarrow (DD(M,P))_\h = DC_{N_g} \coprod DC_{N_g},$$
where $\mathcal{N}_i V$ is the neighborhood of radius $i$ around $V$ with respect to the metric $\rho_i$.  Moreover $DD(M,P)$ has a characteristic submanifold disjoint from $V$.

{\noindent}(2)  One can choose a basepoint in each component of $V$ and $(DD(M,P))_\h$ such that the Riemannian
manifolds $$(V, \rho_i) \quad \text{converge to} \quad (DD(M,P))_\h \stackrel{isom.}{=} DC_{N_g} \coprod
DC_{N_g}$$ in the pointed Gromov-Hausdorff topology.

{\noindent}(3)  $\lim \V (DD(M,P), \rho_i) = \lim \V (V, \rho_i) = 4 \,\V(C_{N_g}).$
\end{thm}

{\noindent}The above theorem says the closed manifold $DD(M,P)$ admits non-positively curved metrics which, in the hyperbolizable pieces of $DD(M,P)$, look more and more like hyperbolic manifolds.  Moreover, the bits of $DD(M,P)$ which are not hyperbolizable become very long and thin, with small total volume.

We will use a version of the barycenter method due to Souto.
\begin{thm} \label{Souto 2'} \cite[App.A,Prop.5]{BCS1}
There exists a sequence of $\mathcal{C}^1$ homotopy equivalences $F_i : A^\infty_i \longrightarrow (DD(M,P),
\rho_i)$ such that: 

{\noindent}(1)  There is a sequence $c_i \searrow 1$ such that $| \text{Jac} \, F_i | \le
c_i$ on $F_i^{-1}(V)$.

{\noindent}(2)  There exist constants $\varepsilon, R, r,$ and $L$ such that: if
$x \in \U{A^\infty_i}$ is such that the ball $B (F_i(x), R) \subset (DD(M,P), \rho_i)$ is hyperbolic, and $|
\text{Jac} \, F_i (x) | \ge (1- \varepsilon) \cdot c_i$, then $F_i$ is $L$-Lipschitz on the ball
$B_{A^\infty_i}(x, r)$.
\end{thm}
{\noindent}Roughly speaking, the maps $F_i$ are nearly volume non-increasing on the bits of $A_i^\infty$ which map into hyperbolic pieces of $(DD(M,P), \rho_i)$.  Moreover, if the map does not dramatically shrink volume at some hyperbolic point in the target, then the map is uniformly Lipschitz.  These two properties can be played off each other to show that the maps $F_i$ are uniformly Lipschitz on compact sets.  The following argument is similar to \cite[Sec.3.1]{BCS1} and \cite[Ch.7]{BCGlong}.

\begin{lem}
There exists an $i_0$ such that for all $i \ge i_0$ the map
\begin{equation} \label{Lipschitz map}
F_i \circ \psi_i \circ \phi^j_i : DC_{N,1/i}^{\ge \zeta} \longrightarrow (DD(M,P), \rho_i)
\end{equation}
is $2L$-Lipschitz, where $j \in \{1,2 \}$.
\end{lem}

\begin{pf}
Without a loss of generality, we may assume that $r < 1$ in Theorem \ref{Souto 2'}.  We have established the limits
$$\V (A^\infty_i) \longrightarrow 4 \text{Vol}(C_N) \quad \text{and} \quad
  \V (DD(M,P), \rho_i) \longrightarrow 4 \V (C_{N_g}).$$ 
The map $F_i$ is nearly volume nonincreasing on the set $F_i^{-1}(V)$.  From this and conclusion (3) of Theorem \ref{Leeb's result} it follows that,
 $$\left( \V(A^\infty_i) - \V (F_i^{-1} (V)) \right) \longrightarrow 0.$$
If $F_i$ uniformly shrunk a large piece of $A_i^\infty$, then too much volume would be lost to hit all of $(V,\rho_i)$ with the nearly volume nonincreasing map $F_i$.  In other words,
$$\V \left( \{ x \in A^\infty_i \ | \ \ | \text{Jac} \, F_i (x) | < (1- \varepsilon) \} \right) 
    \longrightarrow 0.$$
So we've shown that the set of points where $F_i$ shrinks volume by a definite amount is small, and the set of points where we do not control the behavior of $F_i$ is small.  However, these small volume sets may be spread all through the space $A_i^\infty$.  This is where Lemma \ref{Ahlfors reg 2} comes in: for $x \in DC_{N, 1/i}^{\ge \zeta}$, a metric $r$-ball about the point $\psi_i \circ \phi^k_i (x) \in A_i^\infty$ has volume at least $ar^3>0$.  So for $i$ sufficiently large, every point $\psi_i \circ \phi^k_i (x) \in A_i^\infty$ is at most distance $r$ from a point in the (good) set
$$F_i^{-1}(V) \bigcap \left\{ x \in A_i^\infty \ | \ \ | \text{Jac} \, F_i (x) | \ge (1 - \varepsilon ) \right\}.$$
Therefore for $i \gg 0$ an $r$-neighborhood of the set 
$$F_i^{-1} (V) \bigcap \{ x \in A^\infty_i \ | \ \ | 
      \text{Jac} \, F_i (x) | \ge (1- \varepsilon) \}$$ 
contains the subspace 
$$\left( \psi_i \circ \phi_i^1 (DC_{N, 1/i}^{\ge \zeta}) \right) \bigcup
  \left( \psi_i \circ \phi_i^2 (DC_{N, 1/i}^{\ge \zeta}) \right) \subseteq A^\infty_i.$$
($\zeta>0$ is the fixed constant from the beginning of the section.)  From conclusion (2) of Theorem \ref{Souto 2'} we may conclude that for $i \gg 0$, the maps
\begin{equation} 
F_i \circ \psi_i \circ \phi^j_i : DC_{N,1/i}^{\ge \zeta} \longrightarrow (DD(M,P), \rho_i)
\end{equation}
are $2L$-Lipschitz, where $j \in \{1,2 \}$.  This proves the lemma.
\end{pf}

With these geometric facts, we can now use a simple topological argument to prove that $(M,P)$ is pared acylindrical.

Suppose $(M,P)$ is not pared acylindrical.  Let $B \subset M$ be an essential annulus.  Pick a homotopically nontrivial
connected closed curve $\gamma \subset M$ which intersects $B$ essentially.  Recall the homeomorphism 
$$g: (M,P) \longrightarrow (C_N \cap N^{\ge \zeta}, C_N \cap \partial N^{\ge \zeta}).$$ 
Without a loss of generality, we may
assume that $g(\gamma)$ is a smooth finite length curve in $C_N$.  Now consider $\gamma$ as a curve in $DD(M,P)$.  Since
$\gamma \subset M$ intersects $B \subset M$ essentially, $\gamma \subset DD(M,P)$ intersects the torus $DB \subset
DD(M,P)$ essentially.  The complement $\overline{DD(M,P) - V} \subset DD(M,P)$ contains a characteristic
submanifold of $DD(M,P)$, implying $DB$ is homotopic into $\overline{DD(M,P) - V}$.  Therefore $\gamma$ intersects
both $\overline{DD(M,P) - V}$ and $V$ essentially.  This implies the $\rho_i$-length of the geodesic homotopic to
$\gamma$ is going to infinity.  This contradicts the fact that the map (\ref{Lipschitz map}) is $2L$-Lipschitz and
$g(\gamma)$ has finite length.  Therefore $(M,P)$ must be pared acylindrical.   \vskip 6pt

To complete the proof of Theorem \ref{equality case thm}, it remains to prove that $N$ must be minimally
parabolic.  Suppose that $N \in \mathcal{I}(M,P)$ is not minimally parabolic.  As in the proof of Theorem \ref{lower bound thm2}, by enlarging the paring locus $R \supset P$, $N$ becomes minimally parabolic in $\mathcal{I}(M,R)$ (see Section \ref{deformation theory section}).  Since by hypothesis $\V (C_N) = \frac{1}{2} \text{SimpVol}(D(M,P)),$ Theorem \ref{lower bound thm2} implies that $\text{SimpVol}(D(M,P)) = \text{SimpVol}(D(M,R))$.  Therefore the above argument applied to $N$ in
$\mathcal{I}(M,R)$ proves that $(M,R)$ is pared acylindrical.  Thus $D(M,R)$ is hyperbolizable.  But Proposition
\ref{simpvol goes up with drilling 2} then implies that $\text{SimpVol}(D(M,R)) > \text{SimpVol}(D(M,P))$.  This
is a contradiction.  Therefore, $N$ must be minimally parabolic.
\square \vskip 6pt

We can now prove

\vskip 6pt \noindent\textbf{Theorem \ref{BCGthm}.} \itshape  Let $M$ be a compact oriented $3$-manifold with nonempty incompressible boundary such that the interior of $M$
admits a complete convex cocompact hyperbolic metric.  If $N$ is a hyperbolic $3$-manifold homotopy
equivalent to $M$ then $$\text{Vol}(C_N) \ge \frac{1}{2} \text{SimpVol}(DM).$$ Moreover, if $\text{Vol}(C_N) =
\frac{1}{2} \text{SimpVol}(DM) > 0$ then $M$ is acylindrical, $N$ is convex cocompact, and $\partial C_N$ is
totally geodesic. \normalfont \vskip 6pt

\begin{pf}
The inequality follows from Theorem \ref{lower bound thm2}.  So suppose $$\text{Vol}(C_N) = \frac{1}{2}
\text{SimpVol}(DM) > 0.$$ Then by Theorem \ref{equality case thm}, $M$ is acylindrical and $N$ is convex
cocompact.  It follows from Theorem \ref{my thesis} that $\partial C_N$ is totally geodesic.
\end{pf}

\section{Geometric convergence of convex cores} \label{gccc section}

This section will assemble facts from the literature concerning geometric convergence of convex cores.  These facts 
will be used in Section \ref{combinations section}.

We begin with a proposition which follows quickly from the following fact: Kleinian limit sets move continuously in the geometric topology when the injectivity radius at the basepoint is uniformly bounded from above and below \cite{KT}\cite[Prop.2.4]{Mc3}.

\begin{lem} \label{McM fact} Let $(N_n, \omega_n)$ be a sequence of framed hyperbolic $3$-manifolds converging 
geometrically to $(N, \omega)$.  Assume all framed basepoints $\omega_n$ (resp $\omega$) are contained in the 
intersection of the convex core $C_{N_n}$ (resp $C_N$) and the $\mu_3$-thick part $N_n^{\ge \mu_3}$ (resp. $N^{\ge 
\mu_3}$).  Suppose the convex cores $C_{N_n}$ are all homeomorphic, $C_{N_n}$ has incompressible boundary, and there 
is a uniform upper bound on the volume of $C_{N_n}$.  Then the limit sets $L_n$ of $(N_n, \omega_n)$ converge in the 
Hausdorff topology to the limit set $L$ of $(N, \omega)$. \end{lem}

{\noindent}(The framed manifold $(N_n, \omega_n)$ canonically determines a subgroup $\pi_1 (N_n, \omega_n) \subset \text{PSl}_2 \mathbb{C}$.  This notation is defined in Section \ref{deformation theory section}.  The above hypotheses are not optimal.  We prove only what we will need later.)

\begin{pf} 
Using \cite[Prop.2.4]{Mc3}, it suffices to prove there is an $R>0$ such that for any $n$ and any point $x 
\in C_{N_n}$, the injectivity radius of $N_n$ at $x$ is less than $R$.  The uniform upper bound on the volume of the 
convex cores $C_{N_n}$ yields a uniform upper bound $r_1$ on the radius of an embedded ball contained in $C_{N_n}$.  
Since the manifolds $N_n$ are all homeomorphic, it follows from the Gauss-Bonnet theorem that there is a uniform upper bound $r_2$ on the injectivity radii of the surfaces $\partial C_{N_n} \subset N_n$.  Using the fact that $C_{N_n}$ has incompressible boundary, it follows that the injectivity radius at a point $x \in C_{N_n}$ is less than $r_1 + r_2$. \end{pf}

We combine this with the work of Bowditch \cite{Bow} and Taylor \cite[Thm.7.2]{Ta} to obtain

\begin{prop} \label{Taylor fact} Assume the hypotheses of Proposition \ref{McM fact}.  Then the framed submanifolds 
$(C_{N_n} , \omega_n)$ converge in the pointed Gromov-Hausdorff topology to $(C_N , \omega)$. \end{prop}

\begin{pf}  By Proposition \ref{McM fact}, the limit sets $L_n$ of $(N_n, \omega_n)$ converge in the Hausdorff topoogy 
to the limit set $L$ of $(N, \omega)$.  By \cite[Thm.7.2]{Ta}, it follows that for any $r>0$ the submanifolds 
$(C_{N_n} \cap B(\omega_n, r))$ converge in the (basepoint free) Gromov-Hausdorff topology to $(C_N \cap B(\omega, 
r))$.  (In the notation of \cite[Thm.7.2]{Ta}, we know that $\Lambda$ is equal to the full limit set of the geometric 
limit $\Gamma$.  With this, the desired convergence is proven in the second paragraph of the proof of 
\cite[Thm.7.2]{Ta}, which in turn uses a result from \cite{Bow}.)  \end{pf}

We now generalize the geometric topology slightly to handle sequences of manifolds which are pulling apart into 
multiple components.

\begin{defn} \label{geom conv} \emph{(Geometric convergence)} Let $( N_n , \{ \omega^i_n \}_{i=1}^k )$ be a sequence 
of hyperbolic $3$-manifolds equipped with $k$ framed basepoints $\{ \omega^1_n, \ldots \omega^k_n \}$ contained in the 
intersection $C_{N_n} \cap N^{\ge \mu_3}$.  $(N_n , \{ \omega^i_n \} )$ converges geometrically to a hyperbolic 
$3$-manifold $(X, \{ \omega^i \})$ (similarly equipped with $k$ framed basepoints in $C_X \cap X^{\ge \mu_3}$) if: 

{\noindent}$\bullet \ \ $For each compact subset $K \subset X$ containing $\{ \omega^1, \ldots, \omega^k \}$ there are 
smooth \emph{embeddings} $$\psi_n: (K, \omega^1, \omega^2, \ldots, \omega^k) \longrightarrow (N_n, \omega^1_n, \ldots, 
\omega^k_n)$$ converging $\mathcal{C}^\infty$ to an isometric embedding. 

{\noindent}$\bullet \ \ $For $i \neq j$ the distance $d_{N_n} (\omega^i_n, \omega^j_n)$ goes to infinity with $n$.
\end{defn}

Let $(N_n, \{ \omega^i_n \} )$ be a sequence of hyperbolic $3$-manifolds equipped with 
$k$ framed basepoints contained in the intersection $C_{N_n} \cap N^{\ge \mu_3}$.  Assume the sequence $(N_n, \{ 
\omega^i_n \})$ is converging geometrically to $(N, \{ \omega^i \} )$.  
For any $r>0$ define $$K_r := \bigcup_i B_N (\omega^i, r).$$  Then geometric convergence provides a sequence of smooth 
embeddings $$\psi_n: ( K_r ,\omega^1, \omega^2, \ldots, \omega^k) \longrightarrow (N_n, \omega^1_n, \ldots, 
\omega^k_n)$$ converging $\mathcal{C}^\infty$ to an isometric embedding.  

Let us additionally assume that the convex cores $C_{N_n}$ are all 
homeomorphic, $C_{N_n}$ has incompressible boundary, and there is a uniform upper bound on the volume of $C_{N_n}$.  
Then by Proposition \ref{Taylor fact}, for any $\delta>0$ there is an index $n_\delta$ such that $n>n_\delta$ implies: 

{\noindent}$\bullet \ \ $For any $x \in C_N \cap K_r$, $\psi_n (x)$ lies within distance $\delta$ from $C_{N_n}$. 

{\noindent}$\bullet \ \ $For any $y \in C_{N_n} \bigcap \left( \cup_i B_{N_n} (\omega_n^i, r) \right)$, $y$ lies 
within distance $\delta$ from a point in the image $\psi_n (C_N \cap K_r)$. 

{\noindent}In particular, it follows that 
$$\text{Vol}\left( C_{N_n} \bigcap \left( \cup_i B_{N_n} (\omega_n^i, r) \right) \right) \longrightarrow \text{Vol}(C_N \cap K_r),$$ 
implying that 
$$\liminf \text{Vol} (C_{N_n}) \ge \text{Vol}(C_N).$$ 

For later reference let us record this as

\begin{cor} \label{final fact}
Let $(N_n, \{ \omega^i_n \} )$ be a sequence of hyperbolic $3$-manifolds equipped with 
$k$ framed basepoints contained in the intersection $C_{N_n} \cap N^{\ge \mu_3}$.  Assume the sequence $(N_n, \{ 
\omega^i_n \})$ is converging geometrically to $(N, \{ \omega^i \} )$, the convex cores $C_{N_n}$ are all 
homeomorphic, $C_{N_n}$ has incompressible boundary, and there is a uniform upper bound on the volume of $C_{N_n}$.  Then
$$\liminf \text{Vol} (C_{N_n}) \ge \text{Vol}(C_N).$$ 
\end{cor}

\section{Flattening the bending lamination}\label{combinations section}

Recall that $(M,P)$ is a connected pared $3$-manifold with pared incompressible boundary.  To avoid trivial cases,
we assume $\pi_1 (M)$ is not virtually abelian.  The goal of this section is to prove the following theorem.  (For definitions see Section \ref{bending lamination section}.)

\begin{thm} \label{combinations thm 2} Let $\lambda$ be the bending measured lamination on
$\partial M - P$ of the boundary of the convex core of a minimally parabolic geometrically finite hyperbolic manifold $N \in \mathcal{I}(M,P)$ homeomorphic to the interior of $M$.  Then for $0 < \varepsilon < 1$ there exists a geometrically finite minimally parabolic hyperbolic manifold $N_\varepsilon \in \mathcal{I}(M,P)$ homeomorphic to the interior of $M$ with bending measured lamination $\varepsilon \lambda$ such that 
$$\V (C_{N_\varepsilon}) \longrightarrow \frac{1}{2} \text{SimpVol}(D(M,P)) \quad \text{as} 
     \quad \varepsilon \rightarrow 0. \qquad (\dagger)$$
\end{thm}

{\noindent}The existence of the manifolds $N_\varepsilon$ above follows quickly from the work of Bonahon-Otal \cite{BO} when $P = \emptyset$, and Lecuire \cite{Lc} when $P \neq \emptyset$.  (Note the assumption that $(M,P)$ has pared incompressible boundary does not imply that $M$ has incompressible boundary.  Indeed, $(M,P)$ may be a handlebody pared along a multicurve which intersects every compressing disk.  This is why it may be necessary to use results from \cite{Lc}.)  The central point of Theorem \ref{combinations thm 2} is $(\dagger)$.  At first it may seem necessary to pass to a subsequence to obtain $(\dagger)$.  However, we will quickly establish that the volumes $\V (C_{N_\varepsilon} )$ are monotonic in $\varepsilon$, allowing us to pass to subsequences freely in the proof.  

Intuitively, the manifolds $N_{\varepsilon}$ should converge in some sense to a (possibly disconnected) manifold 
with Fuchsian ends, because all the geometric limits of the components of $\partial C_{N_\varepsilon}$ will be totally geodesic surfaces.  Using a fair amount of Kleinian theory, the following argument makes this idea precise, and yields sufficient convergence to establish the limit $(\dagger)$.

\vskip 6pt

The rest of the section is a proof of Theorem \ref{combinations thm 2}.  Pick a measured geodesic lamination
$\lambda$ on $\partial M - P$ realized by a minimally parabolic geometrically finite hyperbolic manifold $N_1$
equipped with a homeomorphism $\phi_1: M-P \longrightarrow C_{N_1}$.

By Theorem \ref{BO thm}, for $0 < \varepsilon < 1$ there exists a geometrically finite hyperbolic manifold
$N_\varepsilon$ such that there is a homeomorphism $\phi_\varepsilon: M - P \longrightarrow C_{N_\varepsilon}$,
and the bending lamination on $\partial C_{N_\varepsilon}$ pulls back via $\phi_\varepsilon$ to $\varepsilon
\lambda$ on $\partial M - P$.  Moreover, the volumes $\V (C_{N_\varepsilon})$ are non-increasing as $\varepsilon \rightarrow 0$.  In particular, by Corollary \ref{final fact} any geometric limit of the manifolds $N_\varepsilon$ (in the sense of Definition \ref{geom conv}) is geometrically finite.  This monotonicity also means that in proving Theorem \ref{combinations thm 2} we are free to pass to subsequences.

We now apply a key compactness result due to Canary-Minsky-Taylor \cite[Thm.5.5]{CMT}.  Let $\rho_\varepsilon : \pi_1 (M) \longrightarrow \text{PSl}_2 \mathbb{C}$ be a representation given by $(N_\varepsilon, \phi_\varepsilon)$ 
($\rho_\varepsilon$ is unique up to conjugation in $\text{PSl}_2 \mathbb{C}$).  In the current notation,
 their theorem implies the following:

\begin{thm} \label{CMTthm} \cite[Thm.5.5]{CMT} There exists a sequence $\varepsilon_n \rightarrow 0$, a sequence of 
homeomorphic re-markings $r_n : (M,P) \longrightarrow (M,P)$ which are the identity on the complement of the window of 
$(M,P)$, a collection $x$ of disjoint, non-parallel, homotopically non-trivial simple closed curves in the window base 
such that:

{\noindent}$(1) \ \ $If $M'$ is a component of $M-X$ whose closure is not a thickened torus (where $X$ is the total 
space of the interval bundle over $x$), then there is a sequence of inner-automorphisms $\sigma_n$ of 
$\text{PSl}_2\mathbb{C}$ such that the sequence of representations $$\{ (\sigma_n \circ \rho_n \circ r_{n*})|_{\pi_1 
(M')} \}$$ converges pointwise.

{\noindent}$(2) \ \ $If $c$ is a curve in $x$ which is not freely homotopic into $P$, then the length of the closed 
geodesic in $N_{\varepsilon_n}$ freely homotopic to $(\phi_n \circ r_n) (c)$ converges to zero.

\end{thm}

{\noindent}For simplicity, let $N_n : = N_{\varepsilon_n}$, $\phi_n := \phi_{\varepsilon_n}$, and $m_n := \phi_n \circ 
r_n$.

Let $( W, R ) \subset M$ be the union of the connected components of the characteristic submanifold of $(M,P)$ which are interval bundles over a surface with a non-abelian fundamental group.  Recall that $(M_{\text{acyl}}, Q)$ is the union of the pared acylindrical connected components of the complement of the characteristic submanifold of $(M,P)$.  It will be convenient to move the components of $X$ into a standard form, and possibly make $X$ larger without violating the statement of Theorem \ref{CMTthm}.  First, if a component of $X$ is isotopic rel $\partial M - P$ into $R \cup Q$, then isotope that component into $R \cup Q$.  (This corresponds to having a component of $x$ which is isotopic into the boundary of the window base.  We want to push it to the boundary of the window base.)  Next, if $A$ is any component of $R \cup Q$, $A$ is not contained in $X$, and $A$ is freely homotopic in $M$ into $X$, then enlarge $X$ by adding the essential annulus $A$.  Finally, if $A$ is any component of $R \cup Q$, $A$ is not contained in $X$, and $A$ is freely homotopic in $M$ into $P$, then enlarge $X$ by adding the essential annulus $A$.  The statement of Theorem \ref{CMTthm} remains true for the enlarged interval bundle $X$.

\subsection{The interval bundle case}

Pick a nonelementary component $M'$ of the $3$-manifold $M - X$.  Let $\U{N}_n$ be the isometric covering of $N_n$ 
corresponding to the subgroup $m_{n*} (\pi_1 (M')) < \pi_1 (N_n)$.  Let $\U{m}_n : M'\longrightarrow \U{N}_n$ denote a 
resulting homotopy equivalence making the diagram

\begin{equation*} \begin{CD} M'   @>{\U{m}_n}>>        \U{N}_n      \\ @VVV        @VVV      \\ M   @>{m_n}>>         
N_n \end{CD} \end{equation*}

{\noindent}commute up to homotopy.  

By conclusion (1) of Theorem \ref{CMTthm}, the sequence $(\U{N_n}, \U{m_n})$ converges in $\text{AH}(M')$.  Let 
$(\U{N},\U{m}) \in \text{AH}(M')$ be the algebraic limit.  To be explicit, pick a basepoint $* \in M' \subset M$, pick 
a framed basepoint $\U{\omega}$ (resp. $\U{\omega}_n$) for $\U{N}$ (resp. $\U{N}_n$), and homotopically alter the 
markings $m_n, \U{m}_n,$ and $\U{m}$ to make them maps of pointed spaces.  By algebraic convergence, these choices can 
be made such that the holonomy representations 

$$\pi_1 (M', *) \stackrel{\U{m}_n}{\longrightarrow} \pi_1 (\U{N}_n, \U{\omega}_n) 
    \stackrel{\text{hol.}}{\longrightarrow} \text{PSl}_2 \mathbb{C}$$

{\noindent}converge pointwise to the holonomy representation

$$\pi_1 (M', *) \stackrel{\U{m}}{\longrightarrow} \pi_1 (\U{N}, \U{\omega}) 
    \stackrel{\text{hol.}}{\longrightarrow} \text{PSl}_2 \mathbb{C}.$$
    
{\noindent}(The choice of a framed basepoint specifies the holonomy representation uniquely.  See Section \ref{deformation theory section}.)  After possibly passing to a subsequence, we may assume that all framed basepoints are contained both in the convex core and in the $\mu_3$-thick part.  Finally, let $\iota_n : \U{N}_n \longrightarrow N_n$ denote the isometric covering map, and define the framed basepoint $\omega_n := \iota_{n *} \U{\omega}_n$.

\begin{prop} \label{I-bundle lem} 
Suppose $(M,P)$ has an essential annulus or M{\"o}bius band which is contained in $M'$ and is \emph{not} 
freely homotopic into $X \cup P$. Then the sequence of framed manifolds $(N_n, \omega_n)$ converges geometrically to the algebraic limit $(\U{N}, \U{\omega})$, and $\U{N}$ is either a Fuchsian or an extended Fuchsian group. \end{prop}

{\noindent}Before proving this lemma, we state two corollaries.

\begin{cor} \label{I-bundle cor}
Suppose $M'$ essentially intersects $W \subseteq M$.  Then $M' \subseteq W$.
\end{cor}

{\noindent}The proof of Corollary \ref{I-bundle cor} is immediate.

\begin{cor} \label{acyl cor}
Suppose $M'$ essentially intersects $M_{\text{acyl}} \subseteq M$.  Then $M'$ is a component of $M_{\text{acyl}}$.
\end{cor}

\emph{Proof of Cor.\ref{acyl cor}:}
Suppose $M'$ essentially intersects $M_{\text{acyl}}$.  Because $(M_{\text{acyl}},Q)$ is pared acylindrical, it follows that a component $(L, \Pi)$ of $(M_{\text{acyl}}, Q)$ is contained in $M'$.

Pick an essential annulus $A$ on the frontier of $L$ as a submanifold of $M'$.  In particular, such an $A$ is not freely homotopic into $X$.  If no such $A$ exists then $M' = L$ and we are done.  If $A$ is freely homotopic (in $M$) into $P$, then by our definition of $X$, $A$ is a component of $X$.  Therefore $A$ is not freely homotopic into $X \cup P$.  We may apply Proposition \ref{I-bundle lem} to conclude that $(L, \Pi)$ is homotopic into the characteristic submanifold of $(M,P)$.  This is a contradiction.  Therefore there is no such essential annulus $A$, implying $M' = L$.
\square \vskip 6pt

\emph{Proof of Prop.\ref{I-bundle lem}:}
The proof begins with a separate claim.  Let $\mathbf{bl}(N_n)$ indicate the bending lamination on the boundary of
the convex core of $N_n$.  Recall that by definition $\mathbf{bl}(N_n) = \phi_{n *} ( \varepsilon_n \lambda)$.
\vskip 4pt

{\noindent}\emph{Claim: For each $n$ there is an essential annulus or M{\"o}bius band $B_n$ of $C_{N_n}$ such that
$m_n^{-1} (B_n) \subset (M,P)$ essentially intersects $M'\subset M$.  Moreover, the intersection number $i(\partial 
B_n, \mathbf{bl}(N_n))$ goes to zero as $n$ goes to infinity.} \vskip 2pt

{\noindent}(The annulus $m_n^{-1}(B_n) \subset M$ may change to follow the re-markings $r_n$.) \\
{\noindent}\emph{Proof of Claim:}  Let $\{ \Sigma_\ell \}_\ell$ be a finite collection of essential annuli in the window of $(M,P)$ which are not freely homotopic into $P$, such that the window minus $\cup_\ell
\Sigma_\ell$ is a collection of topological balls and solid tori, and each component of the window not freely homotopic into $P$ contains a component of the collection $\{ \Sigma_\ell \}_\ell$.  By hypothesis the collection $\{ \Sigma_\ell \}_\ell$ is not empty.  For each pared homeomorphism $r_n :(M,P) \longrightarrow (M,P)$, the submanifold $M'$ will essentially 
intersect $r_n^{-1} (\Sigma_{\ell_n}) \subset (M,P)$ for some ${\ell_n}$.

Let $B_n \subset C_{N_n}$ be the essential annulus or M{\"o}bius band $\phi_n (\Sigma_{\ell_n}) \subset N_n$.  The
intersection number $$i(\partial B_n \, , \, \mathbf{bl}(N_n))
    = i(\phi_n^{-1} (\partial B_n) \, ,  \, \varepsilon_n \lambda)
   = \varepsilon_n \cdot i( \Sigma_{\ell_n} \, , \, \lambda)
   \le \varepsilon_n \cdot \max_\ell \left\{ i( \Sigma_\ell \, , \, \lambda) \right\}$$
goes to zero.  Also, $m_n^{-1}(B_n) = r_n^{-1} (\Sigma_{\ell_n})$ essentially intersects $M'$.  This proves
the claim. \vskip 4pt

Let $\gamma_n^*$ be the closed geodesic freely homotopic into $B_n$.  $B_n \subset C_{N_n}$ is homotopic rel boundary to an essential annulus or M{\"o}bius band $B^*_n$ containing $\gamma_n^*$, such that $\partial B^*_n \subset \partial C_{N_n}$ is a (possibly disconnected) closed geodesic in the intrinsic hyperbolic metric on $\partial C_{N_n}$.  By the claim, $\partial B^*_n$ is very lightly bent, meaning that $i(\partial B^*_n, \mathbf{bf}(N_n))$ is small.  After adjusting $B^*_n$ by an isotopy fixing both $\partial B^*_n$ and $\gamma_n^*$, we may assume that the supremum
$$\sup_{y \in B^*_n} d_{N_n} (y, \gamma_n^*)$$
is attained at a point in $\partial B^*_n$.  (Simply straightening the annulus rel $\partial B^*_n$ and $\gamma_n^*$ would suffice.)

Let $\{ g_n^1, \ldots, g_n^k \}$ be a generating set for $\pi_1 (\U{N}_n, \U{\omega}_n)$.  Let $g_n^{i*}$ be the 
\emph{based} geodesic in $N_n$ corresponding to $g_n^i$.  We chose $B_n$ so that $m_n^{-1}(B_n) \subset M$ would intersect $M'$ essentially.  Therefore one of the based geodesics 
$g_n^{i*}$ of the generating set must intersect $B^*_n$.  By algebraic convergence, the lengths of the based geodesics $\{ g_n^{i*} 
\}_{n,i}$ have a finite supremum.

We now show that the length of $\gamma_n^*$ has a uniform lower bound.  Suppose that up to extracting a subsequence the length
of $\gamma_n^*$ is going to zero.  Since $\partial B^*_n$ is very lightly bent, Lemma \ref{Lecuire lemma} implies that the length of $\partial B^*_n$ is also going to zero.  This forces the entire surface $B^*_n\subset N_n$ to be contained deep inside the Margulis tube corresponding to $\gamma_n^*$.  In particular, this produces a sequence of based geodesics $g_n^{i_n *}$ with unbounded length, contradicting algebraic
convergence.  Therefore the length of $\gamma_n^*$ has a uniform lower bound.

Once we know that $\gamma_n^*$ is uniformly long, the inequality of Lemma \ref{Lecuire lemma} becomes stronger, implying that $\partial B^*_n$ lies inside a $d_n$-neighborhood of $\gamma_n^*$, where $d_n \rightarrow 0$.  By the construction of $B^*_n$ this forces all of $B^*_n$ into a $d_n$-neighborhood of $\gamma_n^*$.  Since the based geodesic $g_n^{i_n *}$ intersects $B^*_n$, it follows that $\gamma_n^* \subset N_n$ passes within a uniformly bounded distance from the basepoint $\omega_n$ of $N_n$.

In the universal cover $\mathbb{H}^3$ with basepoint $o \in \mathbb{H}^3$, we can lift $\gamma_n^*$ to a geodesic passing within a uniformly bounded distance from $o$.  This specifies a lift of the entire surface $B^*_n$ to an infinite closed annulus.  Choose a point on each component of this infinite closed annulus, call them $y^1_n$ and $y^2_n$.  At $y^i_n$ pick a hyperbolic plane $\mathfrak{h}^i_n$ intersecting the convex hull of $\pi_1 (N_n, \omega_n) \subset \mathbb{H}^3$ only at the point $y^i$ (this notation is defined in Section \ref{deformation theory section}).  Such a plane is called a support plane at $y^i$.  We chose $y^1_n$ and $y^2_n$ to lie on different components of the convex hull of $\pi_1 (N_n, \omega_n)$.  This forces the planes $\mathfrak{h}^1_n$ and $\mathfrak{h}^2_n$ to be disjoint.  Moreover, they bound a closed region of $\mathbb{H}^3$ containing the convex hull of $\pi_1 (\U{N}_n, \U{\omega}_n)$ and the convex hull of $\pi_1 (N_n, \omega_n)$.  By choosing the points $y^1_n$ and $y^2_n$ to be as near to $o$ as possible, we can guarantee there is a geodesic segment $\sigma_n$ of length less than $4d_n$ from $\mathfrak{h}^1_n$ to $\mathfrak{h}^2_n$.

Since the geodesic $\gamma_n^*$ passes within a uniformly bounded distance from the basepoint $\omega^1$ of $N_n$, the geodesic segment $\sigma_n$ stretching between the support planes $\mathfrak{h}^1_n$ and $\mathfrak{h}^2_n$ can be 
chosen to stay within a uniformly bounded distance from the basepoint $o \in \mathbb{H}^3$.  Since the length of 
$\sigma_n$ is going to zero, the support planes $\mathfrak{h}^1_n$ and $\mathfrak{h}^2_n$ are converging to each 
other.  Therefore after passsing to a subsequence we can find a round circle at infinity $\mathbb{S}^1 \subset \partial \mathbb{H}^3$ such that: for $k=1,2$ the round circles $\partial \mathfrak{h}^k_n \subset \partial \mathbb{H}^3$ are converging in the Hausdorff topology to $\mathbb{S}^1$.   Therefore the limit sets of the groups $\pi_1 (N_n, \omega_n)$ and $\pi_1 (\U{N}_n, \U{\omega}_n)$ are contained in smaller and smaller neighborhoods of $\mathbb{S}^1 \subset \partial 
\mathbb{H}^3$.  

Suppose that the fundamental group of the algebraic limit $\pi_1 (\U{N}, \U{\omega}) < \text{PSl}_2 \mathbb{C}$ 
contains an element $\gamma$ which does not preserve $\mathbb{S}^1 \subset \partial \mathbb{H}^3$.  Let $U \subset 
\partial \mathbb{H}^3$ be an open neighborhood containing $\mathbb{S}^1$ such that $\gamma.\theta \notin \overline{U}$ for some $\theta \in \mathbb{S}^1$.  We can find a small neighborhood $V \subset U$ of $\theta$ such that for 
sufficiently large $n$, the group $\pi_1 (\U{N}_n, \U{\omega}_n)$ contains an element $\gamma_n$ satisfying $\gamma_n 
. V \cap U = \emptyset$.  Moreover, for sufficiently large $n$ the limit set of $\pi_1 (N_n, \omega_n)$ is contained 
in $U$ and intersects $V$ in at least one point $p$.  (This point $p$ exists because $(M,P)$ has pared incompressible 
boundary, forcing each component of the boundary of the convex hull of $\pi_1 (N_n, \omega_n)$ to be a pleated 
hyperbolic plane.  At least one of these pleated hyperbolic planes limits onto a Jordan curve in $\partial 
\mathbb{H}^3$ which is homotopically nontrivial in the annulus $\overline{\text{Hull}(\mathfrak{h}^1_n \cup 
\mathfrak{h}^2_n)} \cap \partial \mathbb{H}^3$.  This annulus is converging in the Hausdorff topology to 
$\mathbb{S}^1$.)  Since the action of $\pi_1 (N_n, \omega_n)$ on $\partial \mathbb{H}^3$ preserves its own limit set, 
this is a contradiction.  Therefore the fundamental group of the algebraic limit $\pi_1 (\U{N}, \U{\omega})$ preserves 
$\mathbb{S}^1 \subset \partial \mathbb{H}^3$.

Let $\mathbb{H}^2 \subset \mathbb{H}^3$ be the hyperbolic plane spanning $\mathbb{S}^1$.  To see that the algebraic 
limit $\pi_1 (\U{N}, \U{\omega})$ acts on $\mathbb{H}^2$ with finite co-area, consider any component $(S, \partial S)$ 
of $\overline{\partial M - P} \cap \overline{M'}$.  Each component of $\partial S$ is mapped to curves of smaller and 
smaller length (by conclusion (2) of Theorem \ref{CMTthm}).  So in the algebraic limit the surface $(S, \partial 
S)$ must correspond to a quasi-Fuchsian subgroup of $\pi_1 (\U{N}, \U{\omega})$.  Since the full algebraic limit 
$(\U{N}, \U{\omega})$ preserves $\mathbb{S}^1$, and it contains a quasi-Fuchsian subgroup, it must be either Fuchsian 
or extended Fuchsian (with finite co-area).

After possibly passing to a subsequence, let $(Z, \tau)$ be the geometric limit of the framed manifolds $(N_n, \omega_n)$.  By the previous paragraph and Proposition \ref{McM fact}, the Kleinian group $\pi_1 (Z, \tau) < \text{PSl}_2 \mathbb{C}$ corresponding to $(Z, \tau)$ is geometrically finite and preserves the hyperbolic plane $\mathbb{H}^2 \subset \mathbb{H}^3$ spanning $\mathbb{S}^1$.  Since the algebraic limit $\U{N}$ is geometrically finite also, by considering convex core area one sees that $(\U{N}, \U{\omega})$ is a finite cover of the geometric limit $(Z, \tau)$.  The next paragraph will use this to show that the covering map is a homeomorphism.

The locally isometric covering map $(\U{N}, \U{\omega}) \longrightarrow (Z, \tau)$ induces an inclusion $$\pi_1 
(\U{N}, \U{\omega}) \le \pi_1 (Z, \tau) < \text{PSl}_2 \mathbb{C}.$$  Suppose there is an element $g \in \pi_1 (Z, 
\tau ) - \pi_1 (\U{N}, \U{\omega})$.  Since the covering map has finite degree, there is a $k \in \mathbb{N}$ such 
that $g^k \in \pi_1 (\U{N}, \U{\omega})$.  Therefore there exists a sequence of natural numbers $\{ n_i \}$, a 
sequence $\alpha_{n_i}$ in $\pi_1 (M, *)$, and a $\beta \in \pi_1 (M',*)$ such that $$g =\lim_i m_{{n_i}*} 
(\alpha_{n_i}), \quad \text{and} \quad g^k = \lim_i \U{m}_{{n_i}*} (\beta),$$ where $$m_{{n_i}*}: \pi_1 (M,*) 
\longrightarrow \pi_1 (N_n, \omega_n)  \quad \text{and} \quad \U{m}_{{n_i}*}: \pi_1 (M',*) \longrightarrow \pi_1 
(\U{N}_n, \U{\omega}_n)$$ are the maps on fundamental groups induced by the markings.  This implies that $m_{n_i *} 
(\alpha_{n_i}^k \beta^{-1})$ converges to the identity matrix as $i \rightarrow \infty$.  The algebraic convergence of 
the representations $\U{m}_{n*} \rightarrow \U{m}$ forces the existence of an $\eta>0$ such that every basepoint 
$\omega_n \in N_n$ is contained in the $\eta$-thick part of $N_n$.  (Failure would force every loxodromic element of 
$\pi_1 (\U{N}_n, \U{\omega}_n) < \pi_1 (N_n, \omega_n)$ to translate the base point $\U{\omega}_n$ 
arbitrarily far for $n \gg 0$, violating algebraic convergence.)  This uniform thickness $\eta$ yields a uniform 
neighborhood $U \subset \text{PSl}_2 \mathbb{C}$ of the identity such that $$\pi_1 (N_n, \omega_n) \cap U = 
\text{identity} \quad \text{for all  } n.$$  Therefore for $i \gg 0$, $\alpha_{n_i}^k = \beta \in \pi_1 (M,*)$.  In 
the fundamental group of a hyperbolic $3$-manifold, group elements with at least one $k^{\text{th}}$ root have a 
unique $k^{\text{th}}$ root.  So there is an element $\gamma \in \pi_1 (M,*)$ satisfying $\gamma^k = \beta \in \pi_1 
(M',*)$, and $\alpha_{n_i} = \gamma$ for all $i \gg 0$.  It follows that in fact $\gamma \in \pi_1 (M',*)$.  So finally we 
have shown that $g = \lim_i m_{n_i *} (\gamma) = \lim_i \U{m}_{n_i *} (\gamma)$, implying that $g$ is an element of 
the algebraic limit $\pi_1 (\U{N}, \U{\omega})$.  (It is a standard fact from the theory Kleinian groups that if the 
algebraic limit is a finite index subgroup of the geometric limit, then they are the same.  See, for example, the proof of \cite[Prop.4.2]{JM}.  The above argument is an 
adaptation of the proof of this fact to the current setting.)

Thus we have shown that $\pi_1 (Z, \tau) = \pi_1 (\U{N}, \U{\omega}) < \text{PSl}_2 \mathbb{C}$. \square \vskip 6pt

Using Corollaries \ref{I-bundle cor} and \ref{acyl cor} we may conclude that: \vskip 4pt \itshape

{\noindent}Each non-elementary component $M'$ of $M - X$ is either a sub-interval bundle of $( W, R)$, or is equal to a component of $(M_{\text{acyl}}, Q)$.  Conversely, each component of $(M_{\text{acyl}}, Q)$ is a component of $M-X$. \normalfont \vskip 4pt

\subsection{The pared acylindrical case} \label{the pared acylindrical case section}

Pick a component $(L, \Pi)$ from the pared acylindrical manifold $(M_{\text{acyl}},Q)$.  To avoid too much clutter, we will now re-use some previous notation.  We now let $\U{N}_n$ denote an isometric 
covering of $N_n$ corresponding to the conjugacy class of the subgroup $m_{n*} (\pi_1 (L)) < \pi_1 (N_n)$.  Let $\U{m}_n : L \longrightarrow \U{N}_n$ denote a resulting homotopy equivalence making the diagram

\begin{equation*} \begin{CD} L   @>{\U{m}_n}>>        \U{N}_n      \\ @VVV        @VVV      \\ M   @>{m_n}>>         
N_n \end{CD} \end{equation*}

{\noindent}commute up to homotopy.  By Theorem \ref{CMTthm}, the sequence $\{ (\U{N}_n, \U{m}_n) \} \subset 
\text{H}(L)$ converges algebraically to a manifold $(\U{N}, \U{m}) \in \text{H}(L,\Pi)$.

Pick a basepoint $* \in L \subset M$, pick a framed basepoint $\U{\omega}$ 
(resp. $\U{\omega}_n$) for $\U{N}$ (resp. $\U{N}_n$), and homotopically alter the markings $m_n, \U{m}_n,$ and $\U{m}$ 
to make them maps of pointed spaces.  By algebraic convergence, these choices can be made such that the holonomy 
representations 

$$\pi_1 (L, *) \stackrel{\U{m}_n}{\longrightarrow} \pi_1 (\U{N}_n, \U{\omega}_n) 
    \stackrel{\text{hol.}}{\longrightarrow} \text{PSl}_2 \mathbb{C}$$

{\noindent}converge pointwise to the holonomy representation

$$\pi_1 (L, *) \stackrel{\U{m}}{\longrightarrow} \pi_1 (\U{N}, \U{\omega}) 
    \stackrel{\text{hol.}}{\longrightarrow} \text{PSl}_2 \mathbb{C}.$$
    
{\noindent}Up to passing to a subsequence, we may assume that all framed basepoints are contained both in the convex 
core and in the $\mu_3$-thick part.  Finally, let $\iota_n : \U{N}_n \longrightarrow N_n$ denote the isometric 
covering map, and define the framed basepoint $\omega_n := \iota_{n *} \U{\omega}_n$. 

After possibly passing to a subsequence, let $(Z, \tau)$ be a geometric limit of a subsequence of $\{ (N_n, \omega_n ) \}$.  By Proposition \ref{final 
fact}, $\text{Vol}(C_Z) \le \liminf \text{Vol}(C_{ \U{N}_n}) < \infty$.  Therefore $(Z, \tau)$ is geometrically 
finite.  Since there is a locally isometric covering map $(\U{N}, \U{\omega}) \rightarrow (Z, \tau)$, a standard covering theorem due to Thurston now implies that the algebraic limit $\U{N}$ is also 
geometrically finite \cite[Thm.2.1]{CoverC}.

\begin{lem} The algebraic limit of the sequence $(\U{N}_n, \U{m}_n) \in \text{H}(L)$ is minimally parabolic in the set 
$\text{H}(L,\Pi)$. \end{lem} 

\begin{pf}  $(L,\Pi)$ is pared acylindrical.  For some $\zeta < \mu_3 $ we have the homotopy equivalence 
$\U{m}: L \longrightarrow  C_{\U{N}} \cap (\U{N})^{\ge \zeta}$ 
given by the marking.  Using conclusion (2) of Theorem \ref{CMTthm}, we know that (up to homotopy) $\U{m}$ takes $\Pi$ into $C_{\U{N}} \cap \partial (\U{N})^{\ge \zeta}$.  Using a topological rigidity theorem due to Johannson, we can conclude that $\U{m}$
is homotopic to a homeomorphism taking $\Pi$ into $C_{\U{N}} \cap \partial (\U{N})^{\ge \zeta}$ 
\cite[Lem.X.23,pg.235]{Joh}.  Let $\U{m}$ now denote this homeomorphism.

Suppose $(\U{N},\U{m})$ is not minimally parabolic in $\text{H}(L,\Pi)$.  Then there is an essential curve $c_1 \subset 
L$ not freely homotopic into $\Pi$, such that $c_1$ is mapped by $\U{m}$ to a parabolic curve in $\U{N}$. Since $\U{m}$ 
is a homeomorphism, the curve $c_1$ is freely homotopic into $\partial L-\Pi$.

Therefore we may assume there exists an essential curve $c_1 \subset \partial L-\Pi$ such that the length of the closed 
geodesic $\U{m}_n (c_1)^*$ freely homotopic to $\U{m}_n (c_1) \subset \U{N}_n$ is going to zero.  Pick an essential 
curve $c_2 \subset \partial L-\Pi$ which intersects $c_1$ essentially in $\partial L-\Pi$.  Since the re-markings $r_n 
: (M,P) \longrightarrow (M,P)$ are the identity on the pared acylindrical component $(L,\Pi)$, we may apply 
Lemma \ref{Lecuire lemma} to conclude that (for $i=1,2$) $\U{m}_n (c_i)$ is homotopic to a curve in $\partial C_{\U{N}_n}$ whose length is uniformly close (i.e. with small multiplicative and small additive error) to the length of the closed geodesic $\U{m}_n 
(c_i)^*$ homotopic to $\U{m}_n (c_i)$.  Therefore $\U{m}_n (c_1)$ is homotopic to a very short curve in $\partial 
C_{\U{N}_n}$, implying that every curve in $\partial C_{\U{N}_n}$ homotopic to $\U{m}_n (c_2)$ is very long.  
Therefore $\U{m}_n (c_2)^*$ is very long, with length going to infinity.  This contradicts algebraic convergence.  
Therefore there exists no such curve $c_1 \subset \partial L$, and we have proven that the algebraic limit of the 
sequence $(\U{N}_n, \U{m}_n) \in \text{H}(L)$ is minimally parabolic in the set $\text{H}(L,\Pi)$. \end{pf}

\begin{lem} \label{the previous lemma} There is a sequence $a_n \rightarrow \infty$ with the following property: if $c 
\subset N_n$ is a homotopically nontrivial closed path such that $m_n^{-1} (c) \subset M$ essentially intersects $L \subset M$, and is not freely homotopic into $L$, then the length of $c \subset N_n$ is at 
least $a_n$. \vskip 2pt \end{lem} 

\begin{pf} Pick an essential annulus $A \subset M$ on the frontier of $L \subset M$.  By previous 
considerations, $A$ must be a component of $X$.  As a simplifying assumption, let us assume that $A$ is not freely homotopic into $P$.  The proof in the general case requires only slight modifications, and is left to the reader.  Let $(A^{*}_n, \partial A^{*}_n) \subset (C_{N_n}, \partial C_{N_n})$ 
be an embedded essential annulus freely homotopic to $m_n (A) \subset N_n$ such that $A^{*}_n$ contains the closed 
geodesic $\gamma^*_n \subset C_{N_n}$ freely homotopic into $m_n (A)$, and $\partial A^*_n \subset \partial C_{N_n}$ 
is a pair of closed geodesics in the intrinsic metric on $\partial C_{N_n}$.  By pulling tight a homotopy from 
$\partial A^*_n$ to $\gamma^*_n$ we may assume that $A^*_n$ is foliated by curves of length less than the length of 
$\partial A^*_n$.

We claim that the length of $\partial A^*_n$ is going to zero.  To see this, recall that the re-marking homeomorphisms 
$r_n: (M,P) \longrightarrow (M,P)$ are the identity off the window of $(M,P)$.  Therefore off the window, the support 
of the bending lamination $\textbf{bl}(N_n)$ is ``fixed'', i.e. the support of $(m_n^{-1})_* (\textbf{bl} (N_n)) 
\subset \partial M - P$ is constant off the window of $(M,P)$.  Since the total mass of $\textbf{bl}(N_n)$ is going to 
zero, we may again use Lemma \ref{Lecuire lemma} to conclude that the length of each component of $\partial A^*_n$ is close to the length of $\gamma^*_n$.  By conclusion (2) of Theorem 
\ref{CMTthm}, the length of $\gamma^*_n$ is going to zero.  This proves the claim.

For $n \gg 0$, let $T_n \subset N_n$ be the $\mu_3$-Margulis tube about the short closed geodesic $\gamma^*_n \subset 
N_n$.  Since the length of $\partial A^*_n$ is going to zero, $A^*_n$ is foliated by very short curves. Therefore the 
length $a_n$ of the shortest path from $\partial T_n$ to $A^*_n$ goes to infinity.

Now let $c \subset N_n$ be a homotopically nontrivial path such that $m_n^{-1} (c) \subset M$ essentially intersects 
the annulus $A \subset M$.  We claim that the length of $c$ is at least $a_n$.  Without a loss of generality, we 
may assume that $c$ is a closed geodesic in $C_{N_n}$.  Then by topological considerations $c$ cannot be contained in 
$T_n$, and it must intersect $A^*_n$.  

Since the frontier of $L \subset M$ has only a finite number of components, 
performing the above procedure on each of them proves the lemma. \end{pf}

Recall that $(Z, \tau)$ is the geometric limit of the framed 
manifolds $(N_n, \omega_n)$, and $Z$ is geometrically finite by Corollary \ref{final fact}. 

\begin{cor} \label{cor of the previous lemma} 
$(Z, \tau)$ is the geometric limit of the sequence of framed manifolds $(\U{N}_n, \U{\omega}_n)$.
\end{cor}
\begin{pf}
After possibly passing to a subsequence, let $(Z', \tau')$ be the geometric limit of the sequence $(\U{N}_n, \U{\omega}_n)$.  Then there is a locally isometric covering map $(Z', \tau') \rightarrow (Z, \tau)$.  Suppose this map is not a homeomorphism.  Then there is an element in $\pi_1 (Z, \tau)$ which is not in the image of the induced homomorphism on fundamental groups.  This element yields a sequence of elements in $\pi_1 (N_n, \omega_n)$ which are not in the image of the homomorphism $\pi_1 (\U{N}_n, \U{\omega}_n) \rightarrow \pi_1 (N_n, \omega_n)$ induced by the covering map, and are represented by curves $\gamma_n$ based at $\omega_n$ of uniformly bounded length. 

Pick an essential curve $b \subset L$ based at $*$ which is not freely homotopic out of $L$.  Then the closed curve $c_n$ formed by concatenating $b$ and $m_n^{-1} (\gamma_n)$ essentially intersects $L$ and is not freely homotopic into $L$.  Therefore the length of a closed curve in $N_n$ homotopic to $m_n (c_n)$ must have length at least $a_n$ (where $a_n$ is the constant of Lemma \ref{the previous lemma}).  Algebraic convergence provides curves based at $\omega_n$ of uniformly bounded length which are homotopic to $m_n (b)$.  Since $a_n \rightarrow \infty$, this is a contradiction.
\end{pf}

\begin{cor} \label{2nd cor of the previous lemma}
There is a sequence $\ell_n \rightarrow 0$ such that: if $A \subset (M,P)$ is an essential annulus on the frontier of $L$ and $\U{m}_n (A) \subset \U{N}_n$ is freely homotopic into a component $S$ of $\partial C_{\U{N}_n}$ then $\U{m}_n(A)$ is freely homotopic to a curve on $S$ of length less than $\ell_n$.
\end{cor}
{\noindent}Note that because $(L,\Pi)$ is pared acylindrical, $\U{m}_n (A)$ is freely homotopic into either two components of $\partial C_{\U{N}_n}$ (when it is parabolic) or exactly one component of $\partial C_{\U{N}_n}$ (when it is hyperbolic).

\begin{pf}
In the notation of the proof of Lemma \ref{the previous lemma}, $\partial A_n^*$ has length going to zero and lies in $\partial C_{N_n}$.  Lift $\partial A_n^*$ to $\U{N}_n$, where the lifts lie in the complement $\U{N_n} - int(C_{\U{N}_n})$.  Project the lifts down to $\partial C_{\U{N}_n}$ to obtain the desired short curves. \end{pf}

\begin{lem} \label{strong on pared acylindrical pieces} The sequence of framed manifolds $(N_n, \omega_n)$ converges 
geometrically to $(\U{N}, \U{\omega})$. \end{lem}

\begin{pf}  By Corollary \ref{cor of the previous lemma}, it suffices to prove that the sequence $(\U{N}_n, \U{\omega}_n)$ converges geometrically to its algebraic limit $(\U{N}, \U{\omega})$.  Recall that $(\U{N}, \U{m})$ is a geometrically finite minimally parabolic manifold in the deformation set $\text{H}(L,\Pi)$. 

The goal will be to apply a criterion for strong convergence due to Jorgensen-Marden \cite[Sec.4.7]{JM} (see also \cite[Thm.4.2]{Mc1}).  Namely, for each cyclic subgroup $\langle a \rangle < \pi_1 (L,*)$ corresponding to an annular component $A$ of $\Pi$, it suffices to show that the groups $\langle \U{m}_{n*} (a) \rangle < \pi_1 (\U{N}_n, \U{\omega}_n)$ converge geometrically to the rank one cusp $\langle \U{m}_* (a) \rangle < \pi_1 (\U{N}, \U{\omega})$.  If for an infinite number of values $n$ the subgroup $\langle \U{m}_{n*} (a) \rangle$ is a parabolic subgroup, then the conclusion follows from algebraic convergence.  So let us assume without a loss of generality that $\U{m}_{n*}(a)$ is a hyperbolic element of $\pi_1 (\U{N}_n, \U{\omega}_n)$ for all $n$.

As a first step, we must show that a compressible curve on $\partial C_{\U{N}_n}$ cannot be short.  Let $\lambda_n$ be the minimum length of a closed geodesic on a hyperbolic surface which essentially intersects a geodesic of length $\ell_n$ (where $\ell_n$ is the sequence of constants from Corollary \ref{2nd cor of the previous lemma}).  Then $\lambda_n \rightarrow \infty$.  Let $c \subset \partial C_{\U{N}_n}$ be a closed curve which is essential on $\partial C_{\U{N}_n}$ and homotopically trivial in $\U{N}_n$.  Since $(L,\Pi)$ is pared acylindrical, every representative of the free homotopy class $\U{m}_n^{-1}(c) \subset L$ lying on $\partial L$ must intersect $\Pi \subset \partial L$.  Applying Corollary \ref{2nd cor of the previous lemma}, $c$ must intersect an essential curve of $\partial C_{\U{N}_n}$ of length less than $\ell_n$.  Therefore the length of $c$ must be at least $\lambda_n$.  Since $\lambda_n \rightarrow \infty$, this completes the first step.

The isometry $a_n := \U{m}_{n*} (a)$ is hyperbolic, and $(L,\Pi)$ is pared acylindrical.  It follows that $\U{m}_n (A) \subset \U{N}_n$ is freely homotopic into a unique component of $\partial C_{N_n}$.  In fact, $\U{m}_n (A)$ is freely homotopic to a short curve in $\partial C_{\U{N}_n}$ (Corollary \ref{2nd cor of the previous lemma}), and all the compressible curves in $\partial C_{\U{N}_n}$ are long (the previous paragraph).  We may therefore apply a theorem of Bridgeman-Canary \cite[Thm.$1'$]{BrC} to conclude that $a_n$ can be represented in the conformal boundary of $\U{N}_n$ by a short curve, i.e. with length in the Poincar{\'e} metric going to $0$.  


We can now finish the proof using the fairly complete description of geometric convergence (in the geometrically finite setting) given by Proposition 4.7 and Theorem 4.9 of \cite{JM}.  Assume the groups $\langle \U{m}_{n*} (a) \rangle$ converge geometrically to a rank two cusp.  Using the two aforementioned results, the curve $\U{m}_* (a) \subset \U{N}$ is represented in the conformal boundary of $\U{N}$ by a geodesic, i.e. a curve with strictly positive length.  This violates the conclusion of the previous paragraph together with the fact that the domains of discontinuity of the manifolds $\U{N}_n$ converge in the sense of Caratheodory to the domain of discontinuity of $\U{N}$.  (This convergence is conclusion (i) of \cite[Thm.4.9]{JM}.)  We must therefore have the desired strong convergence.

\end{pf}

Thus by examining each nonelementary component of $M - X$ we have proven the following proposition for some $0 < \eta_0 < \mu_3$:
\begin{prop} \label{acyl final prop}
After passing to a subsequence there exists a sequence of sets of $k$ framed basepoints $\{ \omega^i_n \}_{i=1}^k \subset C_{N_n} \cap N^{\ge \eta_0}$, a (possibly disconnected) geometrically finite hyperbolic manifold $N$ with $k$ components, and a set of $k$ framed basepoints $\{ \omega^i \}_{i=1}^k \subset C_N \cap N^{\ge \eta_0}$, such that $(N_n, \{ \omega^i_n \})$ converges geometrically to $(N, \{ \omega^i \} )$ in the sense of 
Definition \ref{geom conv}.  \end{prop}

\subsection{Concluding arguments}

We now recall the notation of Corollary \ref{final fact}.  For any $r>0$ define $$K_r := \bigcup_i B_N (\omega^i, 
r).$$  Then the geometric convergence of Proposition \ref{acyl final prop} provides a sequence of smooth embeddings $$\psi_n: ( K_r ,\omega^1, \omega^2, 
\ldots, \omega^k) \longrightarrow (N_n, \omega^1_n, \ldots, \omega^k_n)$$ converging $\mathcal{C}^\infty$ to an 
isometric embedding.  For any $\delta>0$ there is an index $n_\delta$ such that $n>n_\delta$ implies: 

{\noindent}$\bullet \ \ $For any $x \in C_N \cap K_r$, $\psi_n (x)$ lies within distance $\delta$ from $C_{N_n}$. 

{\noindent}$\bullet \ \ $For any $y \in C_{N_n} \bigcap \left( \cup_i B_{N_n} (\omega_n^i, r) \right)$, $y$ lies 
within distance $\delta$ from a point in the image $\psi_n (C_N \cap K_r)$. 

{\noindent}In particular, it follows that $$\text{Vol}\left(C_{N_n} \bigcap \left( \cup_i B_{N_n} (\omega_n^i, r) 
\right) \right) \longrightarrow \text{Vol}(C_N \cap K_r).$$

\begin{lem}
Each component of $C_N$ is either a totally geodesic surface or a $3$-manifold with totally geodesic boundary.
\end{lem}
\begin{pf}
Pick an essential closed curve $c \subset \partial C_N$ which is a geodesic in the intrinsic hyperbolic metric on $\partial C_N$.  Let $c^* \subset C_N$ be the closed geodesic freely homotopic to $c$.  To prove the lemma it suffices to show that $c$ and $c^*$ have the same length.  Using the almost isometric embeddings $\psi_n$ above, for each manifold $N_n$ we can find a geodesic $c_n^* \subset N_n$ homotopic to $\psi_n (c^*)$.  Homotopic to $c_n^*$ there is a curve $c_n \subset \partial C_{N_n}$ which is a geodesic in the instrinsic metric on $\partial C_{N_n}$.  The inequality of Lemma \ref{Lecuire lemma} implies that $c_n$ and $c_n^*$ are near each other in the Hausdorff metric and
$$\ell(c_n) - \ell(c_n^*) \longrightarrow 0.$$
It follows from the algebraic convergence lemmas above that $\ell (c_n^*) \rightarrow \ell(c^*)$.  By using inverses to the almost isometric embeddings $\psi_n$, it follows that $\liminf \ell(c_n)$ is at least $\ell (c)$.  These facts can be combined to show 
$\ell(c) - \ell (c^*) =0.$
\end{pf}

\begin{lem} \label{final lemma} Pick $\delta>0$.  Define $K_{r,n} := \cup_i B_{N_n} (\omega_n^i, r)$.  There exists $r_0,n_0 >0$ such 
that: if $n > n_0$ and $r> r_0$ then the volume of the set $$C_{N_n} - K_{r,n}$$ is less than $\delta$. \end{lem}

\begin{pf} The main tool of the proof will be the isoperimetric inequality for hyperbolic manifolds with abelian 
fundamental group.  Namely, if $H < \text{PSl}_2 \mathbb{C}$ is abelian and $Y \subset \mathbb{H}^3/H$ is a Lipschitz 
embedded finite volume $3$-dimensional submanifold with boundary then $$\text{Vol}(Y) \le \text{Area}(\partial Y).$$ 
This inequality follows immediately from the proof of \cite[Prop.4.1]{Th3}.

Recall the constant $\eta_0$ from Proposition \ref{acyl final prop}.  Fix $\eta< \eta_0$ sufficiently small such that $C_N \cap N^{> \eta}$ is homeomorphic to $C_N$ via the obvious map which retracts the cusps onto the thick part, the area of the surface of intersection 
$$S : = C_N  \cap \partial N^{\ge \eta}$$ 
is less than $\delta/4$, and the boundary of $S$ has total length less than $\delta/4$.  $S$ is a collection of annuli and tori.  Note that $C_N \cap N^{< \eta}$ has totally geodesic boundary (see Section \ref{deformation theory section}).

Pick $r_0>0$ such that $C_N \cap N^{\ge \eta} \subseteq K_{r_0}.$  Let $r>r_0$.

Using the almost isometric embeddings given by geometric convergence, there is an $n_0$ such that $n > n_0$ 
implies the area of the intersection 
$$S_n := C_{N_n} \cap \partial N_n^{\ge \eta} \cap K_{r,n}$$ 
is less than $\delta/2$, and 
the boundary of $S_n$ has total length less than $\delta/2$.  This estimate can be proven carefully by a packing 
argument, using the fact that the local geometry of $S_n$ is horospherical (i.e. the same for all $n$).  

The manifold $(M,P)$ is ``filled up'' by the components of $N$ in the following sense: each nonelementary interval 
bundle in the characteristic submanifold of $(M,P)$ is formed by gluing up some components of $C_N \cap N^{\ge \eta}$ 
along annuli, and each pared acylindrical component of the complement of the characteristic submanifold of $(M,P)$ is 
homeomorphic to a component of $C_N$.  Since the maps $\psi_n$ are embeddings, for each $n$ the subsurface $\partial 
(C_{N_n} \cap N_n^{\le \eta})$ is formed from $S_n$ by successively performing one of the following two operations: 
joining two boundary components of $S_n$ by a ruled annulus (which lies in the boundary of the convex core $C_{N_n}$), 
or attaching to a single boundary component of $S_n$ a half-infinite ruled annulus (lying on the boundary of $C_{N_n}$) going out a rank one cusp of $N_n$.  

A basic fact from hyperbolic geometry is that the area of a hyperbolic ruled annulus (not necessarily compact) is less 
than the length of its boundary \cite[Thm.9.3.1]{Th}.  Thus the area of $\partial (C_{N_n} \cap N_n^{\le \eta})$ is less than $\delta$ .  So finally, by the isoperimetric inequality, the volume of $(C_{N_n} \cap N_n^{\le \eta})$ is 
at most $\delta$.  \end{pf}

From this it follows that $\text{Vol}(C_{N_n}) \longrightarrow \text{Vol}(C_N)$.  Since each component of $C_N$ is 
either a totally geodesic surface or a $3$-manifold with totally geodesic boundary, it follows from Proposition 
\ref{simpvol 2} that the volume of $C_N$ is exactly half the simplicial volume of the doubled manifold 
$D(\mathcal{N}_1 C_N \cap N^{\ge \eta}, \mathcal{N}_1 C_N \cap \partial N^{\ge \eta})$ (where $\mathcal{N}_1 C_N$ is a radius $1$ regular neighborhood of $C_N$ and $\eta$ is the constant from the proof of Lemma \ref{final lemma}).  Since the manifold $D(\mathcal{N}_1 C_N \cap N^{\ge \eta}, \mathcal{N}_1 C_N \cap \partial N^{\ge \eta})$ can be obtained topologically by cutting $D(M,P)$ along essential tori, it follows from Theorem \ref{simpvol 
theorem} that they have the same simplicial volume.  In particular, it follows that 

$$\text{Vol}(C_{N_n}) \longrightarrow \frac{1}{2} \text{SimpVol} (D(M,P)).$$

{\noindent}This completes the proof of Theorem \ref{combinations thm 2}.  \square \vskip 6pt

\bibliography{post_referee_Storm} \end{document}